\documentclass[a4paper,11pt]{article}
\usepackage{amssymb}
\usepackage{amsmath}
\usepackage{amsfonts}
\usepackage[mathscr]{eucal}
\usepackage{enumerate}
\usepackage{amscd}
\usepackage{indentfirst}
\usepackage{enumerate}
\usepackage{txfonts} 
\newtheorem{thm}{Theorem}[section]
\newtheorem{prop}[thm]{Proposition}

\newtheorem{lem}[thm]{Lemma}
\newtheorem{cor}[thm]{Corollary}
\newtheorem{ex}{Example}[section]
\numberwithin{equation}{section}
\allowdisplaybreaks[4]
\title{DIFFERENTIAL OPERATORS AND CONTRAVARIANT DERIVATIVES IN POISSON GEOMETRY}
\author{Yuji HIROTA}
\date{}
\begin{document}
\maketitle 

\begin{abstract}
We inquire into the relation between the curl operators, the Poisson coboundary operators and 
contravariant derivatives on Poisson manifolds to study the theory of differential operators in Poisson geometry. 
Given an oriented Poisson manifold, we describe locally those two differential operators in terms of Poisson connection whose torsion is vanishing. 
Moreover, we introduce the notion of the modular operator for an oriented Poisson manifold. 
For a symplectic manifold, we describe explicitly the modular operator in terms of the curvature 2-section of Poisson connection, 
analogously to the Weitzenb$\ddot{\rm o}$ck formula in Riemannian geometry. 
\end{abstract}
 
\section{Introduction}

As is well-known in Riemannian geometry~\cite{Bpan03,Jrie02}, the exterior derivative $d$ of the de Rham complex for an oriented Riemannian manifold is written locally as 
\begin{equation}\label{sec1:eqn_exterior derivative}
 d = \sum_i \xi_i\wedge \nabla_{X_i}
\end{equation}
by means of the Levi-Civita connection $\nabla$, a local orthonormal frame field $\{X_i\}$ and its dual coframe field $\{\xi_i\}$. 
One can define another differential operator denoted by $\delta$ as $\delta:=(-1)^\bullet\ast^{-1}\circ d\circ\ast$ by using the Hodge star operator $\ast$. 
Similarly, the operator $\delta$ is also expressed locally by means of $\nabla$ and  $\{X_i\}$ as follows:
\begin{equation}\label{sec1:eqn_codifferential}
 \delta = -\sum_{i}i\,(X_i)\nabla_{X_i}. 
\end{equation}
Composing from $d$ and $\delta$, one can obtain a second order differential operator by $\Delta := d\circ\delta + \delta\circ d$ which is called the Laplace operator.  
$\Delta$ can be described locally in terms of the curvature 2-tensor $\mathcal{R}$ of $\nabla$, which is known as the Bochner-Weitzenb$\ddot{\rm o}$ck formula: 
\begin{equation}\label{sec1:eqn_Weitzenbock}
 \Delta = -\sum_i(\nabla_{X_i}\nabla_{X_i}-\nabla_{\nabla_{X_i}X_i}) - \sum_{i,j}\xi^i\wedge i\,(X_j)\ \mathcal{R}(X_i,\,X_j). 
\end{equation}
In these circumstances, differential operators associated with smooth manifolds might be closely related to connections. 

Let us observe what happens when we apply a similar consideration mentioned above to Poisson manifolds. 
Given a Poisson manifold $(P,\,\Pi)$, one can define a differential operator $\delta_\Pi$ of degree $-1$ 
which will be discussed in Section 2. The complex 
\[
\cdots \overset{\delta_\Pi}{\to} \Omega^{\bullet +1}(P)\overset{\delta_\Pi}{\to}\Omega^\bullet(P)
\overset{\delta_\Pi}{\to}\Omega^{\bullet-1}(P)\overset{\delta_\Pi}{\to}\cdots 
\]
is called the canonical complex~(see \cite{Bdiff88}). It might be interesting to 
consider a map on $\Omega^\bullet(P)$ as $\delta_\Pi\circ d + d\circ\delta_\Pi$ in analogy with the Laplace operator 
in Riemannian geometry. However, the map turns out to be identically zero from easy calculation. Thus we reconsider it by replacing $\Omega^\bullet(P)$ by the 
space of multi-vector fields $\mathfrak{X}^\bullet(P)$. Namely, we suggest a new operator which is the analogue with the Laplace operator by two differential complex 
\[
\cdots \overset{\partial_\Pi}{\to} \mathfrak{X}^{\bullet -1}(P)\overset{\partial_\Pi}{\to}\mathfrak{X}^\bullet(P)
\overset{\partial_\Pi}{\to}\mathfrak{X}^{\bullet+1}(P)\overset{\partial_\Pi}{\to}\cdots 
\]
and 
\[
\cdots \overset{\mathfrak{C}_{\mu}}{\to} \mathfrak{X}^{\bullet +1}(P)\overset{\mathfrak{C}_{\mu}}{\to}\mathfrak{X}^\bullet(P)
\overset{\mathfrak{C}_{\mu}}{\to}\mathfrak{X}^{\bullet-1}(P)\overset{\mathfrak{C}_{\mu}}{\to}\cdots, 
\]
and attempt to describe it in terms of connections for Poisson manifolds and those curvatures, analogously with the Bochner-Weitzenb$\ddot{\rm o}$ck formula. 
In the paper, mainly focusing on the case where a Poisson bivector field is nondegenerate, we describe those differential operators $\partial_\Pi$ and $\mathfrak{C}_{\mu}$ 
locally as (\ref{sec1:eqn_exterior derivative}) and (\ref{sec1:eqn_codifferential}) by means of a contravariant derivative $D$ satisfying some conditions. 
Besides, we introduce a new operator called the modular operator and, in Theorem \ref{sec5:thm_main}, describe it on a symplectic manifold 
by means of 2-tensor induced from $D$ similar to (\ref{sec1:eqn_Weitzenbock}). 

The paper is organized as follows: Section 2 is about some operations frequently utilized in Poisson geometry. 
The above differential operators $\partial_\Pi$ and $\delta_\Pi$ and the Schouten-Nijenhuis bracket are discussed in the section. 
In Section 3, we introduce a new operator, called the modular operator, which depend on the choice of the volume form 
basing on the curl operators for oriented Poisson manifolds and exhibit its examples. 
In addition, it turns out that the modular operator is given by Lie derivative by the modular vector field. 
Section 4 is about the theory of contravariant derivatives for Poisson manifolds. In the section, we show that the Poisson coboundary operator $\partial_\Pi$ for 
Poisson manifold is locally represented in terms of the contravariant derivative. Our main result is exhibited in Section 5. 
We explain in detail the star operator on a symplectic manifold which is studied in \cite{Bdiff88} and show some propositions needed for the proof of the main theorem. 
On the basis of those propositions and results in previous sections, we describe locally the curl operator with respect to the Liouville volume form $\mu$ 
for a symplectic manifold in terms of Poisson connection. Consequently, we obtain an explicit formula of the modular operator with respect to $\mu$, which is analogous to 
the Bochner-Weitzenb$\ddot{\rm o}$ck formula (Theorem \ref{sec5:thm_main}). 
Lastly, we specify the modular operator with respect to a volume form $\mu$ multiplied by a non-vanishing function 
in terms of Poisson connection ~(Theorem \ref{sec5:thm_main2}). 
\vspace{0.5cm}

Throughout the paper, we denote by $\varGamma^\infty(E)$ the space of 
smooth sections of a smooth vector bundle $E$ over a smooth manifold $M$. 
Especially, we use the notation $\mathfrak{X}^k(M)$ and $\Omega^k(M)$ for $\varGamma^\infty(\wedge^kTM)$ and $\varGamma^\infty(\wedge^kT^*M)$, respectively. 
\section{Preliminaries}

We begin with the paper by recalling some fundamental items in Poisson geometry. Let $(P,\Pi)$ be a Poisson manifold. The Poisson bivector field $\Pi$ 
gives rise to a homomorphism 
\begin{equation}\label{sec2:eqn1}
\sharp: \Omega^1(P)\longrightarrow \mathfrak{X}^1(P)
\end{equation}
by $\beta\,(\sharp\alpha)=\langle\beta\wedge\alpha,\,\Pi\rangle$ for any $\alpha,\,\beta \in \Omega^1(P)$. 
So, the Poisson bracket $\{f,\,g\}$ of smooth functions $f,\,g$ is given by 
\[
\{f,g\}:=\langle df\wedge dg,\,\Pi\rangle = (df)(\sharp dg).
\] 
A vector field $\sharp dg$ is called the Hamiltonian vector field of $g$, which is denoted by $H_g$. 
By taking exterior powers of the homomorphism, one can extend it to a $C^\infty(P)$-homomorphism 
\begin{equation}\label{sec2:eqn2}
 \sharp_k:\Omega^k(P)\longrightarrow \mathfrak{X}^k(P),\qquad \sharp_k(\alpha_1\wedge\cdots\wedge\alpha_k) := \sharp\alpha_1\wedge\cdots\wedge\sharp\alpha_k.
\end{equation}
On the other hand, the bivector field $\Pi$ naturally induces a skew-symmetric bilinear map 
\begin{equation}\label{sec2:eqn3}
\Pi_k: \Omega^k(P) \times \Omega^k(P) \longrightarrow C^\infty(P)
\end{equation}
by 
\[
  \Pi_k(\alpha_1\wedge\cdots\wedge\alpha_k,\,\beta_1\wedge\cdots\wedge\beta_k) := \det\,\bigl(\langle\alpha_i\wedge\beta_j,\,\Pi\rangle\bigr)
\]
Those maps (\ref{sec2:eqn2}) and (\ref{sec2:eqn3}) are connected each other by the formula 
\[
\Pi_k(\alpha,\,\beta)=\langle\alpha,\,\sharp_k\beta\rangle,\qquad (\alpha,\,\beta\in \Omega^k(P)) .
\]
When $M$ is a symplectic manifold with a symplectic form $\omega$, the homomorphism (\ref{sec2:eqn1}) is an isomorphism due to the nondegeneracy of $\omega$, 
whose inverse map is a natural map induced from the symplectic form 
\[
 \flat:\mathfrak{X}^1(M)\longrightarrow \Omega^1(M),\qquad \flat\,(X):=i(X)\omega.
\]
Here we denote by $i(X)\omega$ the interior product of $\omega$ by $X\in \mathfrak{X}^1(M)$. 
Similarly to the above, one can extend the map $\flat$ to a $C^\infty(M)$-isomorphism 
\begin{equation*}
 \flat_k: \mathfrak{X}^k(M)\longrightarrow \Omega^k(M),\qquad \flat_k(X_1\wedge\cdots\wedge X_k) := \flat\,(X_1)\wedge\cdots\wedge\flat\,(X_k)
\end{equation*}
and find that it is the inverse map of $\sharp_k$. 
The form $\omega$ induces a natural skew-symmetric bilinear forms $\omega_k$ and $\omega_k^{-1}$ on $\wedge^kTM$ and $\wedge^kT^*M$ respectively: 
\begin{align*}
&\omega_k: \mathfrak{X}^k(M)\times \mathfrak{X}^k(M)\longrightarrow C^\infty(P),\qquad 
\omega_k(X_1\wedge\cdots\wedge X_k,\,Y_1\wedge\cdots\wedge Y_k) := \det\,\bigl(\omega(X_i, Y_j)\bigr),\\
&\omega_k^{-1}: \Omega^k(M)\times \Omega^k(M)\longrightarrow C^\infty(P),\qquad 
\omega_k^{-1}(\alpha,\,\beta) := \omega_k(\sharp_k(\alpha),\,\sharp_k(\beta)).
\end{align*}
We can easily check that $\flat_k$ and $\omega_k$ are connected by 
\begin{equation*}
\omega_k(X,\,Y) = \langle\flat_k(X),\,Y\rangle,\qquad (X,\,Y\in \mathfrak{X}^k(M)). 
\end{equation*}
Therefore, one immediately finds that $\omega_k^{-1}=\Pi_k$ and $\Pi=\sharp_2(\omega)$.
\vspace{0.8cm}

For any multi-vector field on a given smooth manifold $M$, we can define a their product which is called the Schouten bracket \cite{LPVpoi12}. 
If $A\in \mathfrak{X}^a(M),\,B\in \mathfrak{X}^b(M)$, their Schouten bracket $[A, B]_\mathrm{S}$ is given by 
\begin{align*}
&[A,B]_\mathrm{S}(F_1,\,\cdots,\,F_{a+b-1}) \\
 &\qquad := \sum_{\sigma\in \mathfrak{S}_{b,a-1}}\mathrm{sgn}(\sigma)A\ \bigl(B(F_{\sigma(1)},\cdots, F_{\sigma(b)}),\,F_{\sigma(b+1)},\cdots, F_{\sigma(a+b-1)}\bigr)\\
 & \qquad\qquad -(-1)^{(a-1)(b-1)}\sum_{\sigma\in \mathfrak{S}_{a,b-1}}\mathrm{sgn}(\sigma)B\ 
 \bigl(A(F_{\sigma(1)},\cdots, F_{\sigma(a)}),\,F_{\sigma(a+1)},\cdots, F_{\sigma(a+b-1)}\bigr).
\end{align*}
for any $F_1,\cdots, F_{a+b-1}\in C^\infty(M)$. 
Here, the symbol $\mathfrak{S}_{p,q}$ denotes the set of $(p,q)$-shuffles. A $(p,q)$-shuffle $\sigma$ is a permutation of $\{1,2,\cdots,p+q\}$ which 
satisfies $\sigma(1)<\cdots<\sigma(p)$ and $\sigma(p+1)<\cdots<\sigma(p+q)$. 
One can easily check that if $\Pi$ is a Poisson bivector, the Hamiltonian vector field $H_f$ for a smooth function $f$ is represented as 
$H_f=-[f,\,\Pi]_\mathrm{S}$. 

The Schouten bracket defines a graded Lie algebra structure on the spaces of smooth vector fields on $M$ whose degrees are shifted by $-1$. That is, 
it satisfies the following properties: 
for any homogeneous element $A\in \mathfrak{X}^a(M), B\in\mathfrak{X}^b(M)$ and $C\in\mathfrak{X}^c(M)$, 
\begin{enumerate}[\quad(S1)]
\item $[A, B]_{\rm S} = -(-1)^{(a-1)(b-1)}[B, A]_\mathrm{S}$.
\item $(-1)^{(a-1)(c-1)}[A, [B, C]_\mathrm{S}]_\mathrm{S}+(-1)^{(b-1)(a-1)}[B, [C, A]_\mathrm{S}]_\mathrm{S} + (-1)^{(c-1)(b-1)}[C, [A, B]_\mathrm{S}]_\mathrm{S}=0$.
\end{enumerate}
Moreover, the Schouten bracket is shown to be compatible with the wedge product:
\begin{enumerate}[\quad(S1)]\setcounter{enumi}{2}
\item $[A, B\wedge C]_\mathrm{S} = B\wedge [A, C]_\mathrm{S} + (-1)^{(a-1)c}[A, B]_\mathrm{S}\wedge C$.
\item $[A\wedge B, C]_\mathrm{S} = [A, C]_\mathrm{S}\wedge B + (-1)^{a(c-1)}A\wedge [B, C]_\mathrm{S}$.
\end{enumerate}
The Schouten bracket $[\cdot,\,\cdot]_S$ is also characterized by
\begin{equation}\label{sec2:eqn_Schouten}
 i\ ([A,\,B]_\mathrm{S}) = i(A)di(B)-(-1)^{(a-1)(b-1)}i(B)di(A) - (-1)^{(a-1)b}i(A\wedge B)d - (-1)^{a(b-1)}di(A\wedge B)
\end{equation}
in terms of the interior product by a multi-vector field. Here, 
the interior product of $\eta\in \Omega^p(M)$ by $A\in \mathfrak{X}^a(M)$ with $p\geq a$, denoted by $i\ (A)\eta$, is defined to be a differential $(p-a)$-form 
by 
\[
 \langle i\ (A)\eta,\,B\rangle := \langle \eta,\,A\wedge B\rangle 
\]
for any $(p-a)$-vector field $B$. Especially, when $a=0$, $i(A)\eta$ is given by $A\eta$. If $p<a$, then we put $i\ (A)\eta = 0$. 

Similarly, we can define the interior product of a multivector field by a differential form. For any $A\in \mathfrak{X}^a(M)$ and $\eta\in \Omega^p(M)$ with $a\geq p$, 
the interior product of $A$ by $\eta$ is defined to be a $(a-p)$-vector field, denoted by $i\ (\eta)A$, such that 
\[
  \langle \tau,\,i\ (\eta)A\rangle = \langle \tau\wedge\eta,\,A\rangle
\]
for any $(a-p)$-form $\tau$. When $a<p$, we put $i\ (\eta)A=0$. 
\begin{prop}\label{sec2:prop_inner product}
Let $\alpha$ be a differential 1-form on $M$. It holds that 
\[
 i(\alpha)(X\wedge A) = X\wedge (i(\alpha)A) + (-1)^a \alpha(X)A
\]
for any $X\in\mathfrak{X}^1(M)$ and $A\in\mathfrak{X}^a(M)$. 
\end{prop}
{\em Proof}.~ Let $\beta\in \Omega^a(M)$. From the definition, we have 
\begin{align*}
\langle \beta,\,i(\alpha)(X\wedge A)\rangle &= \langle \beta\wedge\alpha,\,X\wedge A\rangle = \langle i(X)(\beta\wedge\alpha),\,A\rangle \\
 &= \langle i(X)\beta \wedge\alpha,\,A\rangle + (-1)^a\langle\beta\wedge i(X)\alpha,\,A\rangle\\
 &= \langle i(X)\beta \wedge\alpha,\,A\rangle + (-1)^a\alpha(X)\langle\beta,\,A\rangle\\ 
 &= \langle \beta,\,X\wedge i(\alpha)A\rangle + (-1)^a\alpha(X)\langle\beta,\,A\rangle.
\end{align*}
Since $\beta$ is any $a$-form, this complete the proof. \qquad\qquad\qquad\qquad\qquad\qquad\qquad\qquad\qquad\qquad\qquad $\Box$
\vspace{0.8cm}

Given a Poisson manifold $(P,\,\Pi)$, one can define homomorphisms $\partial_\Pi^{k}$ from $\mathfrak{X}^k(P)$ to $\mathfrak{X}^{k+1}(P)$, 
called the Poisson coboundary operators, as 
\[
\partial_\Pi^{k}: \mathfrak{X}^k(P)\longrightarrow \mathfrak{X}^{k+1}(P),\qquad \partial_\Pi^{k}(A):=-[A,\,\Pi]_{\rm S}
\]
for each nonnegative integer $k=0,1,2,\cdots$. It is well-known that those operators $\partial_\Pi^{\bullet}$ 
make $(\mathfrak{X}^\bullet(P), \partial^\bullet)$ into a chain complex called the Lichnerowicz complex~(see \cite{DZpoi05}). 
The cohomology of the Lichnerowicz complex is called Poisson cohomology of $P$. 
The operators $\partial_\Pi^k$ and $\sharp_k$ are connected by the following lemma~\cite{DZpoi05,LPVpoi12}: 
\begin{lem}\label{sec2:lemma_Poisson coboundary}
For any $\eta \in\Omega^k(P)$, it holds that 
\[
 \sharp_{k+1}(d\eta) = -[\sharp_k(\eta),\,\Pi]_S = \partial_\Pi^{k}\bigl(\sharp_k(\eta)\bigr). 
\]
\end{lem}

It is known that given a Poisson manifold $(P,\Pi)$, there is another differential complex associated to the Poisson bivector $\Pi$. 
For each $k=1,2,\cdots$, one can define the operator $\delta_\Pi^k:\Omega^k(P)\to \Omega^{k-1}(P)$ as 
\begin{align*}
 &\delta_\Pi^k(\varphi_0\,d\varphi_1\wedge\cdots\wedge d\varphi_k)\notag \\
    &\qquad := \sum_{i=1}^k(-1)^{i+1}\{\varphi_0,\,\varphi_i\}\,d\varphi_1\wedge\cdots \wedge\widehat{d\varphi_i}\wedge\cdots \wedge d\varphi_k \notag \\
   &\qquad\qquad + \sum_{i<j}(-1)^{i+j} \varphi_0d\{\varphi_i,\,\varphi_j\}\wedge d\varphi_1\wedge\cdots\wedge\widehat{d\varphi_i}
 \wedge\cdots\wedge\widehat{d\varphi_j}\wedge\cdots\wedge d\varphi_k, 
\end{align*}
where $\varphi_0,\varphi_1,\cdots,\varphi_k\in C^\infty(P)$ and $\{\cdot,\,\cdot\}$ stands for a Poisson bracket with respect to $\Pi$. 
The operator $\delta_\Pi^\bullet$ was introduced by J.-L. Brylinski in \cite{Bdiff88} and satisfy $\delta_\Pi^{k-1}\circ\delta_\Pi^k=0$. 
The cohomology of the differential complex $(\Omega^\bullet(P),\,\delta_\Pi^\bullet)$ is called the canonical homology or Poisson homology of $P$. 
\section{The Modular Operators}

Let $M$ be an oriented $m$-dimensional smooth manifold and $\mu$ a volume form on $M$. 
For a nonnegative integer $a\,(0\leq k\leq m)$, one can define a $C^\infty(M)$-linear map by
\[
{\mu}^\flat:\mathfrak{X}^a(M)\longrightarrow \Omega^{m-a}(M),
\qquad {\mu}^\flat(A):=i\ (A)\mu. 
\] 
${\mu}^\flat$ is an isomorphism whose inverse map, denoted by ${\mu}^\sharp$, is given by 
$\langle\beta,\,{\mu}^\sharp(\alpha)\rangle :=\langle\beta\wedge \alpha,\,\hat{\mu}\rangle$, where $\hat{\mu}$ stands for the dual vector 
field of $\mu$. Remark that $\langle\mu,\,\hat{\mu}\rangle =1$. 
Consequently, one gets $\mathbb{R}$-linear map $\mathfrak{C}_{\mu}^a$ from $\mathfrak{X}^a(M)$ to $\mathfrak{X}^{a-1}(M)$ by 
\[
 \mathfrak{C}_{\mu}^{a} : \mathfrak{X}^a(M)\longrightarrow \mathfrak{X}^{a-1}(M),\qquad 
 \mathfrak{C}_{\mu}^{a}(A):= ({\mu}^\sharp\!\circ d\!\circ{\mu}^\flat)(A). 
\]
The operator $\mathfrak{C}_{\mu}^{a}$ is called the $a$-th curl operator with respect to $\mu$ on $M$~(see \cite{DZpoi05} and \cite{LPVpoi12}). 
From $d\circ d=0$ it follows that $\mathfrak{C}_{\mu}^{a}\circ\mathfrak{C}_{\mu}^{a+1}=0$. We often use the notation $\mathfrak{C}_{\mu}$ for $\mathfrak{C}_{\mu}^{a}$. 

\begin{ex}
A simple example of the curl operator is the case 
where $M$ is an oriented Riemannian manifold $\mathbb{R}^m$ together with a metric $g$. 
Take the volume form $\mu=\sqrt{\det g}\,dx_1\wedge\cdots\wedge dx_m$ and we have, 
for any vector field $X=\sum_{i=1}^m\xi_i\frac{\partial}{\partial x_i}$ on $\mathbb{R}^m$,
\begin{align*}
 (\mu^\flat\circ\mathfrak{C}_{\mu})(X) &= (d\circ\mu^\flat)(X) 
= d\,\sum_{i=1}^m\xi_i\,i\,\biggl({\frac{\partial}{\partial x_i}}\biggr)(\sqrt{\det g}\,dx_1\wedge\cdots\wedge dx_m)\\
&= \sum_{i=1}^m\ \biggl\{\xi_i\frac{\partial \sqrt{\det g}}{\partial x_i} + \sqrt{\det g}\ \frac{\partial \xi_i}{\partial x_i}\biggr\}\ dx_1\wedge\cdots\wedge dx_m.
\end{align*}
Therefore, the curl operator of $X$ with regard to the volume form $\mu$ is given by 
\[
 \mathfrak{C}_{\mu}(X) = \sum_{i=1}^m\ \biggl\{\xi_i\frac{\partial \sqrt{\det g}}{\partial x_i} + \sqrt{\det g}\ \frac{\partial \xi_i}{\partial x_i}\biggr\}.
\]
In particular, taking the standard metric on $\mathbb{R}^m$, $\mathfrak{C}_{\mu}(X)$ is entirely the divergence of a vector field{\rm :} 
\[
 \mathfrak{C}_{\mu}(X) = \mathrm{div}\ X = \sum_{i=1}^m\ \frac{\partial \xi_i}{\partial x_i}.
\]
\end{ex}

\begin{ex}\label{sec3:ex_symplectic}
We consider a symplectic manifold $\mathbb{R}^2$ with a symplectic form $\omega=\phi(x,y)\ dx\wedge dy$, where $\phi(x,y)\ne 0$ on whole $\mathbb{R}^2$.
Remark that $\omega$ is a volume form on $\mathbb{R}^2$. 
Let $X=f\frac{\partial}{\partial x} + g\frac{\partial}{\partial y}\in \mathfrak{X}^1(\mathbb{R}^2)$. Then, $\mathfrak{C}_{\mu}(X)$ with respect to $\omega$ is given by 
\[
 \mathfrak{C}_{\omega}(X) = \frac{\partial f}{\partial x} + \frac{\partial g}{\partial y} 
- \phi\, \biggl\{f\frac{\partial}{\partial x}\biggl(\frac{1}{\phi}\biggr) +  g\frac{\partial}{\partial y}\biggl(\frac{1}{\phi}\biggr)\biggr\} 
= \mathrm{div}\ X + \frac{1}{\phi}\,X\phi.
\]
\end{ex}

\begin{ex}
Let $(P,\Pi)$ be an oriented $m$-dimensional Poisson manifold with a volume form $\mu$.
A vector field $\Xi_{\mu}\in \mathfrak{X}(P)$ satisfying 
\[
 \mathcal{L}_{\Pi}\mu = i\,(\Xi_{\mu})\mu, \
\]
is called the modular vector field with respect to the volume form $\mu$. Here, $\mathcal{L}_\Pi$ denotes the generalized Lie derivative by $\Pi$ which 
is defined as $\mathcal{L}_\Pi = i(\Pi)\circ d - d\circ i(\Pi)$~(see \cite{LPVpoi12}). 
From the definition of the modular vector field, 
\[
 (\mu^\flat\circ \mathfrak{C}_{\mu})(\Pi) = d(i(\Pi)\mu) = - \mathcal{L}_\Pi\mu = -\mu^\flat(\Xi_{\mu}).
\]
Hence we have 
\begin{equation}\label{sec3:eqn_modular vector fields}
 \mathfrak{C}_{\mu}(\Pi) = -\Xi_{\mu}.
\end{equation}
Since $\mathfrak{C}_{\mu}^{1}\circ \mathfrak{C}_{\mu}^{2}=0$, one finds that the curl operator of the modular vector field is zero~(see Proposition 4.17 in \cite{LPVpoi12}).
\end{ex}

\begin{ex}\label{sec3:ex 3-dimensional case}
 Let $\Pi=f\frac{\partial}{\partial x}\wedge\frac{\partial}{\partial y} + g\frac{\partial}{\partial y}\wedge\frac{\partial}{\partial z}+
 h\frac{\partial}{\partial z}\wedge\frac{\partial}{\partial x}$ be a Poisson bivector on $\mathbb{R}^3$. Namely, $\Pi$ is a smooth bivector field on $\mathbb{R}^3$ such that 
 \begin{equation}\label{sec3:eqn 3-dimensional case}
  f\ \Bigl(\frac{\partial h}{\partial x}-\frac{\partial g}{\partial y}\Bigr) + g\ \Bigl(\frac{\partial f}{\partial y}-\frac{\partial h}{\partial z}\Bigr) 
  + h\ \Bigl(\frac{\partial g}{\partial z}-\frac{\partial f}{\partial x}\Bigr) = 0.
 \end{equation}
 Taking the standard volume form $\mu=dx\wedge dy\wedge dz$, we have 
 \begin{equation}\label{sec3:eqn 3-dimensional case2}
  \mathfrak{C}_{\mu}^{2}(\Pi) = \biggl(\frac{\partial f}{\partial y}-\frac{\partial h}{\partial z}\biggr)\frac{\partial}{\partial x} 
   + \biggl(\frac{\partial g}{\partial z}-\frac{\partial f}{\partial x}\biggr)\frac{\partial}{\partial y} 
   + \biggl(\frac{\partial h}{\partial x}-\frac{\partial g}{\partial y}\biggr)\frac{\partial}{\partial z}.
 \end{equation}
\end{ex}
\vspace{0.8cm}

From the formula (\ref{sec2:eqn_Schouten}), 
\begin{align*}
 i([A,B]_\mathrm{S})\mu &= i(A)di(B)\mu -(-1)^{(a-1)(b-1)} i(B)di(A)\mu + (-1)^{a(b-1)}di(A\wedge B)\mu\\
 &= i(\mathfrak{C}_{\mu}(B)\wedge A)\mu - (-1)^{(a-1)(b-1)}i(\mathfrak{C}_{\mu}(A)\wedge B)\mu + (-1)^{a(b-1)}i(\mathfrak{C}_{\mu}(A\wedge B))\mu 
\end{align*}
for any $A\in\mathfrak{X}^a(M),\,B\in\mathfrak{X}^b(M)$. Therefore, we have 
\begin{equation}\label{sec3:eqn_Schouten and curl}
 (-1)^{(a-1)(b-1)}[A,\,B]_{\mathrm{S}} = (-1)^b\mathfrak{C}_{\mu}(A\wedge B) - \mathfrak{C}_{\mu}(A)\wedge B - (-1)^bA\wedge \mathfrak{C}_{\mu}(B).
\end{equation}
We remark that the sign convention of the above formula (\ref{sec3:eqn_Schouten and curl}) is different from the one which Koszul displayed in \cite{Acrochet85}. 
By using (\ref{sec3:eqn_Schouten and curl}) and $\mathfrak{C}_{\mu}\circ\mathfrak{C}_{\mu}=0$, we obtain the following useful proposition, immediately.
\begin{prop}\label{sec3:prop_Schouten and curl}
Let $\mu$ be a volume form on $M$. For any a-vector field $A$ and any b-vector field $B$,  
\[
 \mathfrak{C}_{\mu}[A,\,B]_{\mathrm{S}} = [\mathfrak{C}_{\mu}(A),\, B]_\mathrm{S} + (-1)^{a-1}[A,\, \mathfrak{C}_{\mu}(B)]_\mathrm{S}.
\]
\end{prop} 
\vspace*{0.8cm}

Let $(P,\,\Pi)$ be an oriented Poisson manifold of dimension $n$ with a volume form $\mu$ on $P$. 
For a nonnegative integer $a~(0\leq a\leq n)$, we define a map $\Lambda_{\mu}^{a}$ from $\mathfrak{X}^a(P)$ to $\mathfrak{X}^a(P)$ 
by composing the differential operator $\partial_\Pi^a$ and the curl operator $\mathfrak{C}_{\mu}^a$ with respect to $\mu$ as follows:
\[
 \Lambda_{\mu}^{a}:\mathfrak{X}^a(P)\longrightarrow \mathfrak{X}^a(P),\qquad \Lambda_{\mu}^{a}A:= 
(\mathfrak{C}_{\mu}^{a+1}\!\circ\partial_\Pi^{a}-\partial_\Pi^{a-1}\!\circ\mathfrak{C}_{\mu}^{a})\,A.
\]
Especially, $\Lambda_{\mu}^{0} f = \mathfrak{C}_{\mu}^{1}(H_f)$ for $f\in C^\infty(P)$. 
For the sake of simplicity, we often omit a superscript $a$ and write $\Lambda_{\mu}$ for $\Lambda_{\mu}^{a}$. 
The operator $\Lambda_{\mu}^a$ is a kind of measures the commutativity of 
$\mathfrak{C}_{\mu}^{a+1}\!\circ\partial_\Pi^{a}$ and $\partial_\Pi^{a-1}\!\circ\mathfrak{C}_{\mu}^{a}$. 

\begin{ex}\label{sec3:ex_2-dimensional case}
A 2-dimensional real space $\mathbb{R}^2$ equipped with a bivector 
$\Pi=\phi(x,y)\frac{\partial}{\partial x}\wedge \frac{\partial}{\partial y}\,(\phi\in C^\infty(\mathbb{R}^2))$ 
is an oriented Poisson manifold with the standard volume form $\mu=dx\wedge dy$. 
The Hamiltonian vector field $H_F$ for $F\in C^\infty(\mathbb{R}^2)$ is expressed as 
\[
H_F = \phi(x,y)\,\biggl(\frac{\partial F}{\partial y}\frac{\partial}{\partial x} 
  - \frac{\partial F}{\partial x}\frac{\partial}{\partial y}\biggr). 
\]
Then, by (\ref{sec3:eqn_modular vector fields}) we have 
\[
 \Lambda_{\mu}^{0} F = \frac{\partial\phi}{\partial x}\frac{\partial F}{\partial y} - 
\frac{\partial\phi}{\partial y}\frac{\partial F}{\partial x} = -\mathfrak{C}_{\mu}^2(\Pi)F = \Xi_{\mu}F. 
\]
and compute $\Lambda_{\mu}^1X$ of a vector field $X=f\frac{\partial}{\partial x}+g\frac{\partial}{\partial y}$ as 
\begin{align*}
 \Lambda_{\mu}^{1}X &= \left(\frac{\partial\phi}{\partial y}\frac{\partial f}{\partial x}- \frac{\partial f}{\partial y}\frac{\partial\phi}{\partial x}
 -f\frac{\partial^2\phi}{\partial y\partial x}-g\frac{\partial^2\phi}{\partial y^2}\right)\frac{\partial}{\partial x} 
 + \left(-\frac{\partial\phi}{\partial x}\frac{\partial g}{\partial y} + \frac{\partial g}{\partial x}\frac{\partial\phi}{\partial y}
 + f\frac{\partial^2\phi}{\partial x^2} + g\frac{\partial^2\phi}{\partial y\partial x}\right)\frac{\partial}{\partial y}\\
 &= -\left(\Xi_{\mu}f+\mathcal{L}_X\frac{\partial\phi}{\partial y}\right)\frac{\partial}{\partial x} 
 + \left(-\Xi_{\mu}g + \mathcal{L}_X\frac{\partial\phi}{\partial x}\right)\frac{\partial}{\partial y}.
\end{align*}
\end{ex}

\begin{ex}\label{sec3:ex_3-dimensional case 2}
Let us consider the Poisson manifold in Example \ref{sec3:ex 3-dimensional case}. The curl operator of $\Pi$ with respect to the standard volume form 
is given by (\ref{sec3:eqn 3-dimensional case2}). From $\partial_\Pi\Pi = 0$ and $(\ref{sec3:eqn 3-dimensional case})$, $\Lambda_{\mu}^{2}\Pi = 0$. 
That is, $A=\Pi$ is a solution of $\Lambda_{\mu}^{2}A=0$. By the way, the Hamiltonian vector field $H_F$ of $F\in C^\infty(P)$ is explicitly 
expressed as 
\[
 H_F = \biggl(f\frac{\partial F}{\partial y}-h\frac{\partial F}{\partial z}\biggr)\frac{\partial}{\partial x}
 + \biggl(g\frac{\partial F}{\partial z}-f\frac{\partial F}{\partial x}\biggr)\frac{\partial}{\partial y}
 + \biggl(h\frac{\partial F}{\partial x}-g\frac{\partial F}{\partial y}\biggr)\frac{\partial}{\partial z}. 
\]
Accordingly, 
\[
\Lambda_{\mu}^0F 
 = \biggl(\frac{\partial h}{\partial z}-\frac{\partial f}{\partial y}\biggr)\frac{\partial}{\partial x} 
   + \biggl(\frac{\partial f}{\partial x}-\frac{\partial g}{\partial z}\biggr)\frac{\partial}{\partial y} 
   + \biggl(\frac{\partial g}{\partial y}-\frac{\partial h}{\partial x}\biggr)\frac{\partial}{\partial z}.
\]
By (\ref{sec3:eqn_modular vector fields}) and (\ref{sec3:eqn 3-dimensional case2}), we find that $\Lambda_{\mu}^0F=\Xi_{\mu}F$. 
\end{ex}

In Example \ref{sec3:ex_2-dimensional case} and \ref{sec3:ex_3-dimensional case 2}, the operator $\Lambda_{\mu}^0$ for a smooth function is described 
in terms of the modular vector field. Those results generally hold for any oriented Poisson manifold: 

\begin{prop}\label{sec3:prop_0-th Poisson-Laplace}
Let $(P,\,\Pi)$ be an oriented Poisson manifold with a volume form $\mu$ and $\Xi_{\mu}$ the modular vector field with respect to $\mu$. 
For any $F\in C^\infty(P)$, it holds that 
\[
 \Lambda_{\mu}^0 F = \Xi_{\mu}F.
\]
\end{prop}
{\em Proof}.~ Put $A=F\in C^\infty(P),\,B=\Pi\in\mathfrak{X}^2(P)$ in (\ref{sec2:eqn_Schouten}) and we have 
\begin{equation*}
 i\,([F,\,\Pi]_\mathrm{S})\mu = Fd i\,(\Pi)\mu - d(i\,(F\Pi)\mu). 
\end{equation*}
Therefore, 
\begin{align*}
 H_F &= ({\mu}^\sharp\circ{\mu}^\flat)(H_F)= -{\mu}^\sharp\bigl(i\,([F,\,\Pi]_\mathrm{S})\mu\bigr)\\ 
     &= {\mu}^\sharp d{\mu}^\flat(F\Pi) - F{\mu}^\sharp d {\mu}^\flat(\Pi)\\
     &= \mathfrak{C}_{\mu}^2(F\Pi) - F\mathfrak{C}_{\mu}^2(\Pi). 
\end{align*}
Since $\mathfrak{C}_{\mu}^1\circ\mathfrak{C}_{\mu}^2=0$, 
\[
  \Lambda_{\mu}^0 F = \mathfrak{C}_{\mu}^1(H_F) = -\mathfrak{C}_{\mu}^1\bigl(F\mathfrak{C}_{\mu}^2(\Pi)\bigr) = -{\mu}^\sharp\bigl(df\wedge i\,(\mathfrak{C}_{\mu}^2(\Pi))\mu\bigr).
\]
Computing the interior product by $\mathfrak{C}_{\mu}^2(\Pi)\in\mathfrak{X}(P)$ of $df\wedge \mu=0$, we find 
\[
 df\wedge i\,(\mathfrak{C}_{\mu}^2(\Pi))\mu = \bigl((\mathfrak{C}_{\mu}^2\Pi)F\bigr)\mu.
\]
From this and (\ref{sec3:eqn_modular vector fields}), it follows that 
\[
 \Lambda_{\mu}^0 F = -\mathfrak{C}_{\mu}^2(\Pi)F\langle\mu,\,\hat{\mu}\rangle = \Xi_{\mu}F, 
\]
which complete the proof. \qquad\qquad\qquad\qquad\qquad\qquad\qquad\qquad\qquad\qquad\qquad\qquad\qquad\qquad\qquad $\Box$

\noindent Proposition \ref{sec3:prop_0-th Poisson-Laplace} gives us a geometric interpretation of $\Lambda_{\mu}^0$: 
$\Lambda_{\mu}^0F$ describe the rate of change of the function along the flow of the modular vector field. 
Moreover, applying Proposition \ref{sec3:prop_Schouten and curl} to the case where $B$ is a Poisson bivector $\Pi$ and using (\ref{sec3:eqn_modular vector fields}), 
we find that 
\begin{equation}\label{sec3:eqn_modular operator}
\Lambda_{\mu}^a A = \mathfrak{C}_{\mu}^{a+1}(\partial^a_\Pi A) - \partial^{a-1}_\Pi(\mathfrak{C}_{\mu}^{a}A) = (-1)^a[\Xi_{\mu},\,A]_\mathrm{S}.
\end{equation}
That is the reason why we call the operator $\Lambda_{\mu}^a$ the $a$-th modular operator with respect to $\mu$. We often omit the superscript $a$ and use the notation 
$\Lambda_{\mu}$ for $\Lambda_{\mu}^a$. 
\section{Poisson Connections}
Let $(P,\Pi)$ be an $n$-dimensional Poisson manifold. As is mentioned in Section 2, the Poisson bivector $\Pi$ induces a homomorphism 
$\sharp:\Omega^1(P)\to \mathfrak{X}^1(P)$. 
By a contravariant derivative on $(P, \Pi)$ we mean an $\mathbb{R}$-linear map 
\[
 D : \Omega^1(P) \longrightarrow \mathfrak{X}^1(P,T^*P):= \varGamma^\infty(TP\otimes T^*P) 
\]
which satisfies the following: for any $f\in C^\infty(P)$ and $\alpha,\eta\in\Omega^1(P)$ 
 \begin{enumerate}[\quad(1)]
  \item $\langle f\eta,\,D\alpha\rangle = f\langle\eta,\,D \alpha\rangle\in \Omega^1(P) ;$
  \item $\langle \eta,\,D(f\alpha)\rangle = f\langle\eta,\,D \alpha\rangle + \bigl((\sharp\eta)f\bigr)\alpha$.
 \end{enumerate} 
Here, $\langle\cdot,\,\cdot\rangle$ stands for the evaluation given by $\langle\eta,\,X\otimes\alpha\rangle:=\eta(X)\alpha$ for any 
$\eta\in\Omega^1(P), X\otimes\alpha\in\mathfrak{X}^1(P,T^*P)$. 
We also often write $D_\eta\alpha$ for $\langle\eta, D \alpha\rangle$ in what follows. 
Similarly to the case of a usual connection theory, $D$ is a local operator. 
Let $\{\theta^1,\cdots,\theta^n\}$ be a local coframe field on a local chart $U$ in $P$. 
If a 1-form $\alpha$ is written on $U$ in the form $\alpha = \sum_i \alpha_i\theta^i~(\alpha_i\in C^\infty(U))$, then 
$D \alpha$ is represented as 
\begin{equation*}
 D \alpha = \sum_{j=1}^n\Bigl\{-\sharp(d\alpha_j) + \sum_{i=1}^n\alpha_i\Theta_{ij}\Bigr\}\otimes \theta^j, 
\end{equation*}
where $\Theta_{ij}$ are local smooth vector fields by $D\theta^i=\sum_{j=1}^n\Theta_{ij}\otimes \theta^j$. 
A contravariant derivative $D$ can be extended to the map from $\mathfrak{X}^a(P)$ to 
$\mathfrak{X}^1(P,\wedge^aTP)= \varGamma^\infty(TP\otimes \wedge^aTP)$ by 
\begin{equation}\label{sec4:contravariant derivative}
 \langle\alpha_1\wedge\cdots\wedge\alpha_a,\,D_\eta A\rangle = \bigl(\sharp\eta\bigr)\bigl(\langle\alpha_1\wedge\cdots\wedge\alpha_a,\,A\rangle\bigr) - 
 \sum_{i=1}^k\langle\alpha_1\wedge\cdots \wedge D_\eta\alpha_i\wedge\cdots\wedge\alpha_a,\,A\rangle
\end{equation}
for any 1-form $\alpha_1,\cdots,\alpha_a$ on $P$. We use the same notation $D$ for the extended contravariant derivative. One can verify that $D$ satisfies 
\begin{equation}\label{sec4:eqn_distribution low}
 D_\eta(A\wedge B) = (D_\eta A)\wedge B + A\wedge (D_\eta B) 
\end{equation}
for any $A\in \mathfrak{X}^a(P),\,B\in \mathfrak{X}^b(P)$ and $\beta\in\Omega^1(P)$. Moreover, by using (\ref{sec4:contravariant derivative}), 
we can get easily the following:
\begin{prop}\label{sec4:prop_inner product}
 Let $\alpha,\,\beta$ be any differential 1-form on $P$. It holds that 
\[
 D_\eta\circ i\ (\alpha) = i\ (\alpha)D_\eta + i\ (D_\eta\alpha). 
\]
\end{prop}

For the Poisson bivector field $\Pi$, $D\Pi$ is given by 
\begin{equation}\label{sec4:contravariant derivative2}
 \langle\beta\wedge\gamma,\,D_\alpha\Pi\rangle  = (\sharp\alpha)\bigl(\langle\beta\wedge\gamma, \Pi\rangle\bigr) - \langle D_\alpha\beta\wedge\gamma,\,\Pi\rangle 
    - \langle\beta\wedge D_\alpha\gamma,\,\Pi\rangle,
\end{equation}
where $\alpha,\beta$ and $\gamma$ is any 1-form on $P$.
The operator $D$ is called a Poisson connection if $D\Pi=0$. In other words, a Poisson connection is the operator $D$ which satisfies
\begin{equation}\label{sec4:Poisson connection}
(\sharp\alpha)\bigl(\langle\beta\wedge\gamma,\,\Pi\rangle\bigr) = \langle D_\alpha\beta\wedge\gamma,\,\Pi\rangle + \langle\beta\wedge D_\alpha\gamma,\,\Pi\rangle.
\end{equation}
It is known that every Poisson manifolds has a Poisson connection (see \cite{Fcon00}). 

As is well-known, one can define a Lie bracket denoted by $[\cdot,\cdot]_\Pi$ on $T^*P$ as 
\[
 [\alpha,\beta]_\Pi:= \mathcal{L}_{\sharp\alpha}\beta - \mathcal{L}_{\sharp\beta}\alpha + d\langle \alpha\wedge\beta,\,\Pi\rangle
\]
for any $\alpha, \beta\in \Omega^1(P)$. Here, $\mathcal{L}_{\sharp\alpha}\beta$ denotes a Lie derivative of $\beta$ in the direction of $\sharp\alpha$. 
One can verify that $[df,\,dg]_\Pi = -d\{f,\,g\}$, and moreover 
the map $\sharp:\Omega^1(P)\to \mathfrak{X}^1(P)$ is a Lie algebra homomorphism: $\sharp[df,\,dg]_\Pi=[\sharp df,\,\sharp dg]$. 
By using this bracket, one can construct a 2-vector field $R_D$ with values in $TP\otimes T^*P$ by putting 
\[
 R_D(\alpha,\beta) := D_\alpha D_\beta - D_\beta D_\alpha - D_{[\alpha,\beta]_\Pi},
\]
which is called the curvature 2-section of $D$. Additional to this, one can get a tensor field $T_D$, called a torsion of $D$, by 
\[
 T_D(\alpha, \beta) := D_\alpha\beta - D_\beta\alpha - [\alpha,\beta]_\Pi. 
\]
The torsion $T_D$ is said to be vanishing if $T_D(\alpha,\beta)=0$ for any $\alpha,\beta\in \Omega^1(P)$. 

Suppose that $(P,\Pi)$ has a Poisson connection $D$ whose torsion is vanishing. 
By the definition of the Schouten-Nijenhuis bracket and (\ref{sec4:contravariant derivative2}), for any $k$-vector field $A$ on $P$ and any function 
$F_1,\cdots,F_{k+1}\in C^\infty(P)$, $(\partial_\Pi^{k} A)(F_1,\cdots,F_{k+1})$ is calculated to be 
\begin{align*}
 &(\partial_\Pi^{a} A)(F_1,\cdots,F_{a+1}) \\
=\,& -[A,\,\Pi]_{\mathrm{S}}(F_1,\cdots,F_{a+1})\\
 =\,& -\!\!\sum_{\sigma\in\mathfrak{S}_{2,a-1}}\!\!\mathrm{sgn}(\sigma) A\bigl(\{F_{\sigma(1)}, F_{\sigma(2)}\},F_{\sigma(3)},\cdots,F_{\sigma(a+1)}\bigr)
     + (-1)^{a-1}\!\!
    \sum_{\sigma\in\mathfrak{S}_{a,1}}\!\!\mathrm{sgn}(\sigma) \bigl\{A(F_{\sigma(1)},\cdots,F_{\sigma(a)}), F_{\sigma(a+1)}\bigr\}\\
 =\,& \sum_{i<j}(-1)^{i+j}A\bigl(\{F_i, F_j\},F_1,\cdots,\widehat{F_i},\cdots,\widehat{F_j},\cdots,F_{a+1}\bigr)
    +\sum_{i=1}^{a+1} (-1)^{i}\bigl\{A(F_1,\cdots,\widehat{F_i},\cdots,F_{a+1}), F_i\bigr\}\\
 =\,& \sum_{i<j}(-1)^{i+j}A\bigl(\{F_i, F_j\},F_1,\cdots,\widehat{F_i},\cdots,\widehat{F_j},\cdots,F_{a+1}\bigr) 
     + \sum_{i=1}^{a+1}(-1)^{i}(D_{dF_i}A)(dF_1,\cdots,\widehat{dF}_i,\cdots,dF_{a+1})\\
   &\qquad + \sum_{i\ne j}^{a+1}(-1)^{i}\langle dF_1\wedge\cdots\wedge\widehat{dF}_i\wedge\cdots\wedge D_{dF_i}dF_j\wedge\cdots\wedge dF_{a+1},\,A\rangle \\ 
 =\,& \sum_{i<j}(-1)^{i}\cdot (-1)^{j-2}
    \langle D_{dF_i}dF_j\wedge dF_1\wedge\cdots\wedge \widehat{dF}_i\wedge\cdots\wedge \widehat{dF}_j\wedge\cdots\wedge dF_{a+1},\,A\rangle \\
    &\qquad +\sum_{i>j}(-1)^{i}\cdot (-1)^{j-1}
    \langle D_{dF_i}dF_j\wedge dF_1\wedge\cdots\wedge \widehat{dF}_j\wedge\cdots\wedge \widehat{dF}_i\wedge\cdots\wedge dF_{a+1},\,A\rangle\\
    &\qquad + \sum_{i<j}(-1)^{i+j}A\bigl(\{F_i, F_j\},F_1,\cdots,\widehat{F_i},\cdots,\widehat{F_j},\cdots,F_{a+1}\bigr)
       + \sum_{i=1}^{a+1}(-1)^{i}(D_{dF_i}A)(dF_1,\cdots,\widehat{dF}_i,\cdots,dF_{a+1})\\
 =\,& \sum_{i<j}(-1)^{i+j}
    \langle D_{dF_i}dF_j\wedge dF_1\wedge\cdots\wedge \widehat{dF}_i\wedge\cdots\wedge \widehat{dF}_j\wedge\cdots\wedge dF_{a+1},\,A\rangle\\
    &\qquad +\sum_{i<j}(-1)^{i+j}
    \langle -D_{dF_j}dF_i\wedge dF_1\wedge\cdots\wedge \widehat{dF}_j\wedge\cdots\wedge \widehat{dF}_i\wedge\cdots\wedge dF_{a+1},\,A\rangle\\
    &\qquad +\sum_{i<j}(-1)^{i+j}A\bigl(\{F_i, F_j\},F_1,\cdots,\widehat{F_i},\cdots,\widehat{F_j},\cdots,F_{a+1}\bigr)
     + \sum_{i=1}^{a+1}(-1)^{i}(D_{dF_i}A)(dF_1,\cdots,\widehat{dF}_i,\cdots,dF_{a+1})\\
 =\,& \sum_{i<j}(-1)^{i+j}
    \langle(D_{dF_i}dF_j-D_{dF_j}dF_i+d\{F_i,F_j\})\wedge dF_1\wedge\cdots\wedge \widehat{dF}_i\wedge\cdots\wedge \widehat{dF}_j\wedge\cdots\wedge dF_{a+1},\,A\rangle\\
    &\qquad +\sum_{i=1}^{a+1}(-1)^{i}(D_{dF_i}A)(dF_1,\cdots,\widehat{dF}_i,\cdots,dF_{a+1})
\end{align*}
  From the assumption, it follows that $T(dF_i,\,dF_j) = D_{dF_i}dF_j-D_{dF_j}dF_i+d\{F_i,F_j\}=0$ for each $i,j=1,\cdots,a+1$. Therefore, it holds that 
 \begin{equation}\label{sec4:eqn Poisson coboundary operator}
  (\partial_\Pi^{a} A)(F_1,\cdots, F_{a+1}) = \sum_{i=1}^{a+1}(-1)^{i} (D_{dF_i}A)(dF_1,\cdots, \widehat{dF}_i,\cdots,dF_{a+1}).
 \end{equation}

\noindent The following proposition is checked easily from (\ref{sec4:eqn Poisson coboundary operator}).
\begin{cor}
$DA=0$ implies that $\partial_\Pi A=0$.
\end{cor}

Let $\{e_1,\cdots,e_n\}$ be a local frame fields on some local chart $U$ of $(P,\Pi)$ and $\{\theta^1,\cdots,\theta^n\}$ be the dual coframe field of 
$\{e_1,\cdots,e_n\}$. Namely,$\{e_1(x),\cdots,e_n(x)\}$ is a basis of $T_xU$ at each point $x\in P$ and $\langle\theta^i(x),\,e_j(x)\rangle=\delta_{ij}~(i,j=1,\cdots,n)$. 
\begin{thm}\label{sec4:thm_Poisson-Lichnerowics}
On some local chart $U$ of $(P,\Pi)$, $\partial_\Pi^{k}$ can be written in the form 
\[
\partial_\Pi^{a}|_U A  = -\sum_{i=1}^{n} e_i\wedge D_{\theta^i}A.
\]
by using a Poisson connection $D$ whose torsion is vanishing. 
\end{thm}
\noindent {\em Proof.}~ Let $A\in \mathfrak{X}^a(P)$ and $F_1,\cdots,F_{a+1}\in C^\infty(P)$. 
\begin{align*}
 &\ \sum_{i=1}^{n}e_i\wedge D_{\theta^i}A(dF_1,\cdots,dF_{a+1})\\ 
&\qquad\qquad\qquad= \frac{1}{a!}\sum_{i=1}^n\sum_{\sigma\in\mathfrak{S}_{a+1}}\mathrm{sgn}(\sigma)(dF_{\sigma(1)})(e_i)
   \bigl(D_{\theta^i}A\bigr)(dF_{\sigma(2)},\cdots,dF_{\sigma(a+1)})\\
&\qquad\qquad\qquad= \frac{1}{a!}\sum_{\sigma\in\mathfrak{S}_{a+1}}\mathrm{sgn}(\sigma)
   \bigl(D_{dF_{\sigma(1)}}A\bigr)(dF_{\sigma(2)},\cdots,dF_{\sigma(a+1)})\\
&\qquad\qquad\qquad= \sum_{i=1}^a(-1)^{i+1} (D_{dF_i}A)(dF_1,\cdots, \widehat{dF}_i,\cdots,dF_{a+1}).
\end{align*}
By the formula (\ref{sec4:eqn Poisson coboundary operator}), we complete the proof. \qquad\qquad\qquad\qquad\qquad\qquad\qquad\qquad\qquad\qquad\qquad $\Box$ 

Theorem \ref{sec4:thm_Poisson-Lichnerowics} does not depend on the choice of the local frame field and its local coframe field. 
Indeed, we let $\{f_1,\cdots,f_n\}$ be another local frame field with the dual coframe field $\{\xi^1,\cdots, \xi^n\}$. 
Then, it holds that 
\[
 f_j = \sum_{i=1}^n A_{ij}e_i\qquad \text{and}\qquad \xi^j = \sum_{k=1}^n B_{jk} \theta^k\qquad (j=1,2,\cdots,n)
\]
for some real coefficients $\{A_{ij}\}$ and $\{B_{jk}\}$. We remark that 
\[
  \sum_{k=1}^n A_{ik}B_{kj} = \delta_{ij} \qquad (1\leq i,j\leq n).
\]
Consequently, 
\[
 \sum_{k=1}^n f_k\wedge D_{\xi^k} = \sum_{k=1}^n\ \biggl(\sum_{i=1}^n A_{ik}e_i\biggr)\wedge D_{\sum_{j=1}^n B_{kj}\theta^j} 
  = \sum_{i,j}\sum_{k=1}^n\ A_{ik}B_{kj} e_i\wedge D_{\theta^j} = \sum_{i=1}^n e_i\wedge D_{\theta^i}.
\]
Hence, $\partial_\Pi$ is written locally as Theorem \ref{sec4:thm_Poisson-Lichnerowics} independently of the choice of the local frame field.

\begin{ex} 
 Let $P$ be an oriented Poisson manifold of $n$-dimension and $A\in\mathfrak{X}^0(P)=C^\infty(P)$.  Then we can check easily 
 \[
  (\partial_\Pi^0 A)(F)= H_AF = -(\sharp dF)A = -D_{dF}A 
 \] 
 for any smooth function $F$. That is, the equation (\ref{sec4:eqn Poisson coboundary operator}) surely holds. Moreover, 
 We choose a local chart $U$ with coordinates $(x_1,\cdots,x_n)$ and a local frame field $\{e_j\}$ by $e_j=\frac{\partial}{\partial x_j}$ on $U$. 
 Then, we have 
 \begin{align*}
  -\sum_{i=1}^{n} e_i\wedge D_{\theta^i}A &= -\sum_{i=1}^{n} \frac{\partial}{\partial x_i}\wedge (\sharp dx_i)A = \sum_{i=1}^n dx_i(H_A)\ \frac{\partial}{\partial x_i} 
   = H_A\\  &= \partial_\Pi A, 
 \end{align*}
 which indicates Theorem \ref{sec4:thm_Poisson-Lichnerowics}. 
\end{ex}
\section{The Modular Operators on Symplectic Manifolds}
\subsection{The star operators for even dimensional Poisson manifolds}

We suppose that $(P,\Pi)$ is an oriented Poisson manifold of dimension $2m~(m\in\mathbb{N})$. As a volume form of $M$, we take a $2m$-form $\mu$ such that 
$\Pi_{2m}(\mu,\,\mu)=1$ and then, we set $\hat{\mu}_P = \sharp_{2m}(\mu)$. 
We consider a natural map 
\[
 \star : \Omega^k(P)\longrightarrow \Omega^{2m-k}(P),\qquad \star \eta := i\,(\sharp_k \eta)\,\mu \qquad\quad (0\leq k\leq 2m)
\] 
which is composed of $\mu^\flat$ and $\sharp_k$. We can verify easily the following:
\begin{equation}\label{sec5:eqn1}
 \star f = f\mu \quad (f\in C^\infty(P)),\qquad\qquad \star\mu = 1.
\end{equation}
\begin{lem}\label{sec5:lemma1}
Let $0\leq k\leq 2m$. If $\alpha\in\Omega^k(P),\,\beta\in\Omega^{2m-k}(P)$, it hold that 
\begin{equation*}
 \Pi_{2m-k}(\star\alpha,\,\beta)\,\mu = \alpha\wedge \beta. 
\end{equation*}
\end{lem}
{\em Proof.}~ Suppose that $\alpha\in\Omega^k(P)$ and $\beta\in\Omega^{2m-k}(P)$ are written in the form respectively 
 \[
   \alpha = \alpha_1\wedge\cdots\wedge\alpha_k,\qquad \beta=\beta_1\wedge\cdots\wedge\beta_{2m-k},
 \]
 where $\alpha_1,\cdots,\alpha_k,\beta_1,\cdots, \beta_{2m-k}$ are differential 1-forms on $P$. 
 Then, we have 
\begin{align}
 \Pi_{2m-k}(\star\alpha,\,\beta) &= \langle\star\alpha,\,\sharp_{2m-k}(\beta)\rangle\notag \\
 &= \mu\bigl(\sharp_k(\alpha),\,\sharp_{2m-k}(\beta)\bigr) = \langle\mu,\,\sharp_k(\alpha)\wedge\sharp_{2m-k}(\beta)\rangle \notag \\
 &= \langle\mu,\, \sharp_{2m}(\alpha_1\wedge\cdots\wedge\alpha_k\wedge\beta_1\wedge\cdots\wedge\beta_{2m-k})\rangle \notag \\
 &= \langle\alpha\wedge\beta,\,\sharp_{2m}(\mu)\rangle \in C^\infty(P). \label{sec5:eqn1 in lemma}
\end{align}
From $\Pi_{2m}(\mu,\,\mu)=1$ it follows that $\Pi_{2m-k}(\star\alpha,\,\beta)\,\mu 
= \langle\alpha\wedge\beta,\,\sharp_{2m}(\mu)\rangle\,\mu = \alpha\wedge\beta$. \qquad\qquad\qquad $\Box$ 

The equation (\ref{sec5:eqn1 in lemma}) leads us to the formula $\star(\alpha\wedge\beta)=\Pi_{2m-k}(\star\alpha,\,\beta)$. 
By using (\ref{sec5:eqn1}) and Lemma \ref{sec5:lemma1}, we calculate 
\[
 \star\star(\alpha\wedge\beta) = \star\Pi(\star\alpha,\,\beta) = \Pi(\star\alpha,\,\beta)\mu = \alpha\wedge\beta. 
\]
Since differential $k$-forms of this type generates $\Omega^k(P)$, we get the following lemma:
\begin{lem}\label{sec5:lemma2}
 For any $\eta\in\Omega^k(P)$, it holds that $\star(\star\eta) = \eta$.
\end{lem}

\begin{prop}\label{sec5:prop1}
 For $\alpha,\,\beta\in\Omega^k(P)$, 
\begin{equation*}
 \alpha \wedge\star\beta = \Pi_k(\alpha,\,\beta)\,\mu = (-1)^k\,\beta\wedge\star\alpha 
\end{equation*}
\end{prop}
{\em Proof}.~ By Lemma \ref{sec5:lemma1} and Lemma \ref{sec5:lemma2}, we calculate 
 \begin{align*}
  \alpha\wedge\star\beta &= (-1)^{k(2m-k)}\star\beta\wedge\alpha = (-1)^{k(2m-k)}\Pi_k(\star(\star\beta),\,\alpha)\,\mu\\
   &= (-1)^{k(2m-k)}\Pi_k(\beta,\,\alpha)\,\mu = \Pi_k(\alpha,\,\beta)\,\mu,
 \end{align*}
 which proves the former equality. The latter follows immediately from $\Pi_k(\alpha,\,\beta)=(-1)^k\Pi_k(\beta,\,\alpha)$. \qquad$\Box$ 

\begin{prop}\label{sec5:prop_relation}
For any $\alpha\in\Omega^k(P)$, it holds that ${\mu}^\sharp(\alpha)=(\sharp_{2m-k}\circ \star)(\alpha)$. 
Especially, if $P$ is a symplectic manifold, it holds that ${\mu}^\flat(K) = (\star\circ\flat_k)(K)$ for any $K\in \mathfrak{X}^k(P)$. 
\end{prop}
{\em Proof}.~ Let $\alpha\in\Omega^k(P)$. Then, by Lemma \ref{sec5:lemma1} we calculate  
\begin{align*}
\langle \beta,\,{\mu}^\sharp(\alpha)\rangle &= \langle\beta\wedge\alpha,\,\hat{\mu}\rangle = (-1)^k\langle\alpha\wedge\beta,\,\hat{\mu}\rangle 
 = (-1)^{2m}\cdot (-1)^{-k}\Pi_{2m-k}(\star\alpha,\,\beta)\\ 
&= \langle \beta,\,\sharp_{2m-k}(\star\alpha)\rangle 
\end{align*}
for any $\beta\in\Omega^{2m-k}(P)$. Therefore, ${\mu}^\sharp = \sharp_\bullet\circ\star$ holds. If $P$ is a symplectic manifold, then the map 
$\sharp_\bullet:\mathfrak{X}^\bullet(P)\to \Omega^\bullet(P)$ is invertible. By Lemma \ref{sec5:lemma2}, the star operator $\star$ is also invertible with 
$\star^{-1}=\star$. Therefore, 
\[
 {\mu}^\flat(K) = ({\mu}^\sharp)^{-1}(K) = (\sharp_k\circ \star)^{-1}(K) = (\star\circ\flat_k)(K),
\]
 which proves the proposition. \qquad\qquad\qquad\qquad\qquad\qquad\qquad\qquad\qquad\qquad\qquad\qquad\qquad\qquad $\Box$

\noindent The following result can be proved immediately by a direct computation:

\begin{prop}\label{sec5:prop_relation2}
 For any $\alpha\in\Omega^1(P)$ and $\eta\in\Omega^k(P)$, it holds that 
\[
 \star(\alpha\wedge\eta) = (-1)^ki\ (\sharp\alpha)\star\eta. 
\]
\end{prop}

In the case where $P$ is a symplectic manifold, the operator $\star:\Omega^\bullet(M)\to\Omega^{2m-\bullet}(P)$ coincides with the one introduced 
by J. -L. Brylinski which is 
defined  as (1) in Proposition \ref{sec5:prop1}. The following proposition is due to J. -L. Brylinski~\cite{Bdiff88}:
\begin{prop}\label{sec5:prop_Brylinski}
 For $0\leq k\leq 2m$, the operator $\delta_\Pi^k:\Omega^k(P)\to \Omega^{k-1}(P)$ is represented as 
 \[
   \delta_\Pi^k = (-1)^{k+1}\,\star\circ d\circ\star. 
 \]
\end{prop}

\subsection{The curl operators on symplectic manifolds } 
We consider a symplectic manifold $(P,\omega)$ of dimension $2m$ with the volume form $\mu=\frac{1}{m!}\wedge^m\omega$ and assume that 
it has a Poisson connection $D$ whose torsion is vanishing. We write $\Pi_\omega$ for a Poisson bivector by the symplectic form $\omega$ and denote by $\sharp_\omega$ 
a map from $\Omega^1(P)$ to $\frak{X}^1(P)$ induced from $\Pi_\omega$~(see Section 2).
Using $D$, we define a covariant derivative $\nabla_XY\in \mathfrak{X}^1(P)$ as 
\begin{equation}\label{sec5:eqn2}
  \nabla_XY := (\sharp_\omega\circ D_{\flat(X)}\circ \flat)(Y)\qquad (X, Y\in \mathfrak{X}^1(P)) 
\end{equation}
and extend it to the map from $\Omega^k(P)$ to $\Omega^1(P, \wedge^kT^*P)$ as 
\begin{equation}\label{sec5:eqn_symplectic connection}
 \nabla_X\eta := (\flat_k\circ D_{\flat(X)}\circ\sharp_k)(\eta)
\end{equation}
for any $\eta\in \Omega^k(P)$. We can verify that the extended map $\nabla_X$ satisfies 
\[
 (\nabla_X\eta)(X_1,\cdots,X_k) = X\,(\eta(X_1,\cdots,X_k)) - \sum_{i=1}^k\eta\,(X_1,\cdots,\nabla_XX_i,\cdots, X_k).
\]
and 
\[
  \nabla_X(\xi\wedge\eta) = \nabla_X\xi\wedge \eta + \xi\wedge\nabla_X\eta\qquad (\xi,\,\eta\in\Omega^\bullet(P)).
\]
For each differential 1-form $\alpha,\,\beta,\,\gamma$ on $P$, there exist vector fields $X,\,Y,\,Z\in \mathfrak{X}^1(P)$ such that 
$\alpha =\flat(X),\,\beta=\flat(Y)$ and $\gamma=\flat(Z)$. Then, we have 
\[
 \sharp_\omega(\alpha)(\Pi_\omega(\beta,\,\gamma)) = X\bigl(\Pi_\omega(\flat(Y),\,\flat(Z))\bigr) = X\bigl(\omega(Y, Z)\bigr).
\]
On the other hand, 
\begin{align*}
 \Pi_\omega(D_\alpha\beta,\,\gamma) + \Pi(\beta,\,D_\alpha\gamma) &= \Pi_\omega(D_{\flat(X)}\flat(Y),\,\flat(Z)) + \Pi_\omega(\flat(Y),\,D_{\flat(X)}\flat(Z))\\
 &= -\flat(Z)\bigl(\sharp_\omega D_{\flat(X)}\flat(Y)\bigr) + \flat(Y)\bigl(\sharp_\omega D_{\flat(X)}\flat(Z)\bigr)\\
 &= \omega (\nabla_XY,\, Z) + \omega (Y,\,\nabla_XZ). 
\end{align*}
Hence, the condition that $D$ is a Poisson connection is equivalent to 
\begin{equation}\label{sec5:eqn3}
 X\bigl(\omega(Y, Z)\bigr) = \omega (\nabla_XY,\, Z) + \omega (Y,\,\nabla_XZ)\qquad (X, Y, Z\in \mathfrak{X}^1(P)).
\end{equation}
The differential 2-form $\omega$ is said to be parallel if (\ref{sec5:eqn3}) holds~(see \cite{Fcon00}). 
On the other hand, the torsion $\nabla$ which is denoted by $T_\nabla$ is vanishing. Indeed, it follows from $T\equiv 0$ that 
\begin{align*}
 T_\nabla(X, Y) &= \nabla_XY - \nabla_YX - [X,\,Y] \\
 &= \sharp_\omega\bigl(D_{\flat(X)}\flat(Y)\bigr) - \sharp_\omega\bigl(D_{\flat(Y)}\flat(Z)\bigr) - [\sharp_\omega(\alpha),\,\sharp_\omega(\beta)]\\
 &= \sharp_\omega(T(\alpha,\,\beta))\\
 &= 0.
\end{align*}
In general, a covariant derivative on a symplectic manifold is called a symplectic connection if the symplectic form is parallel and its torsion is vanishing 
(see \cite{BCGRSsym06}). 
Therefore, we get the following:

\begin{prop}
The covariant derivative $\nabla$ by (\ref{sec5:eqn2}) is a symplectic connection. 
\end{prop}
The condition (\ref{sec5:eqn3}) or (\ref{sec4:Poisson connection}) is also represented as 
\begin{equation}\label{sec5:Poisson connection}
 (D_{\flat(X)}\alpha)(\sharp_\omega\beta) - (D_{\flat(X)}\beta)(\sharp_\omega\alpha) = X\bigl(\Pi_\omega(\alpha,\,\beta)\bigr)
 \qquad (\alpha, \beta\in \Omega^1(P),\,X\in\mathfrak{X}^1(P)).
\end{equation}
Moreover, it can be extended for differential $k$-forms as follows:
\begin{lem}\label{sec5:lemma3}
For any $\alpha, \beta\in \Omega^k(P)$ written in the form $\alpha = \alpha_1\wedge\cdots\wedge \alpha_k,\, \beta=\beta_1\wedge\cdots\wedge\beta_k$ and 
any $X\in \mathfrak{X}^1(P)$, it holds that 
\begin{align*}
 &\sum_{i=1}^k(\alpha_1\wedge\cdots \wedge D_{\flat(X)}\alpha_i\wedge\cdots \wedge\alpha_k)(\sharp_\omega\beta_1,\cdots,\sharp_\omega\beta_k) \\
 &\qquad\qquad + (-1)^k\sum_{j=1}^k(\beta_1\wedge\cdots\wedge D_{\flat(X)}\beta_j\wedge\cdots\wedge \beta_k)(\sharp_\omega\alpha_1,\cdots,\sharp_\omega\alpha_k)
 = X\Bigl(\Pi_\omega(\alpha,\,\beta)\Bigr).
\end{align*}
\end{lem}
{\em Proof}.~ The second term on the left hand side is calculated to 
\begin{align*}
&(-1)^k\sum_{j=1}^k(\beta_1\wedge\cdots\wedge D_{\flat(X)}\beta_j\wedge\cdots\wedge \beta_k)(\sharp_\omega\alpha_1,\cdots,\sharp_\omega\alpha_k)\\
&\quad = 
(-1)^k\sum_{j=1}^k\sum_{\tau\in\mathfrak{S}_k}\mathrm{sgn}(\tau)\ \beta_1(\sharp_\omega\alpha_{\tau(1)})\cdots D_{\flat(X)}\beta_j(\sharp_\omega\alpha_{\tau(j)})
\cdots\beta_k(\sharp_\omega\alpha_{\tau(k)})\\
&\quad = 
-\sum_{j=1}^k\sum_{\tau\in\mathfrak{S}_k}\mathrm{sgn}(\tau)\ \alpha_{\tau(1)}(\sharp_\omega\beta_1)\cdots D_{\flat(X)}\beta_j(\sharp_\omega\alpha_{\tau(j)})
\cdots\alpha_{\tau(k)}(\sharp_\omega\beta_k)\\
&\quad = 
-\sum_{j=1}^k\sum_{\sigma\in\mathfrak{S}_k}\mathrm{sgn}(\sigma)\ \alpha_{\sigma^{-1}(1)}(\sharp_\omega\beta_1)\cdots D_{\flat(X)}\beta_j(\sharp_\omega\alpha_{\sigma^{-1}(j)})
\cdots\alpha_{\sigma^{-1}(k)}(\sharp_\omega\beta_k)\\
&\quad = 
-\sum_{i=1}^k\sum_{\sigma\in\mathfrak{S}_k}\mathrm{sgn}(\sigma)\ \alpha_1(\sharp_\omega\beta_{\sigma(1)})\cdots 
D_{\flat(X)}\beta_{\sigma(i)}(\sharp_\omega\alpha_i)\cdots\alpha_k(\sharp_\omega\beta_{\sigma(k)}). 
\end{align*}
In addition, by using (\ref{sec5:Poisson connection}), we have 
\begin{align*}
 &\sum_{i=1}^k(\alpha_1\wedge\cdots \wedge D_{\flat(X)}\alpha_i\wedge\cdots \wedge\alpha_k)(\sharp_\omega\beta_1,\cdots,\sharp_\omega\beta_k) \\
 &\qquad\qquad + (-1)^k\sum_{j=1}^k(\beta_1\wedge\cdots\wedge D_{\flat(X)}\beta_j\wedge\cdots\wedge \beta_k)(\sharp_\omega\alpha_1,\cdots,\sharp_\omega\alpha_k)\\
=& \sum_{i=1}^k\sum_{\sigma\in\mathfrak{S}_k}\mathrm{sgn}(\sigma)\ 
 \bigl\{D_{\flat(X)}\alpha_i(\sharp_\omega\beta_{\sigma(i)}) - D_{\flat(X)}\beta_{\sigma(i)}(\sharp_\omega\alpha_i)\bigr\}\ 
\alpha_1(\sharp_\omega\beta_{\sigma(1)})
   \cdots \widehat{\alpha_i(\sharp_\omega\beta_{\sigma(i)})}\cdots\alpha_k(\sharp_\omega\beta_{\sigma(k)})\\ 
=& \sum_{i=1}^k\sum_{\sigma\in\mathfrak{S}_k}\mathrm{sgn}(\sigma)\ 
   X\bigl(\Pi_\omega(\alpha_i,\,\beta_{\sigma(i)})\bigr)\Pi_\omega(\alpha_1,\beta_{\sigma(1)})\cdots \widehat{\Pi_\omega(\alpha_i,\beta_{\sigma(i)})}
   \cdots \Pi_\omega(\alpha_k,\beta_{\sigma(k)})\\
=& \sum_{\sigma\in\mathfrak{S}_k}\mathrm{sgn}(\sigma)\ X\biggl(
\Pi_\omega(\alpha_1,\beta_{\sigma(1)})\cdots \Pi_\omega(\alpha_i,\beta_{\sigma(i)})\cdots \Pi_\omega(\alpha_k,\beta_{\sigma(k)})
\biggr)\\
=& X\Bigl(\det\ \bigl(\Pi_\omega(\alpha_i,\beta_j)\bigr)\Bigr),  
\end{align*}
which proves the lemma. 
\qquad\qquad\qquad\qquad\qquad\qquad\qquad\qquad\qquad\qquad\qquad\qquad\qquad\qquad $\Box$
\begin{prop}\label{sec5:prop_Poisson connection}
For any differential $k$-form $\alpha,\,\beta\in\Omega^k(P)$ and any vector field $X\in\mathfrak{X}^1(P)$, it holds that 
\[
\Pi_k(\nabla_X\alpha,\,\beta) + \Pi_k(\alpha,\,\nabla_X\beta) = \nabla_X\bigl(\Pi_k(\alpha,\beta)\bigr).
\]
\end{prop}
{\em Proof}.~ The left hand side is given by 
\[
\Pi_k(\nabla_X\alpha,\,\beta) + \Pi_k(\alpha,\,\nabla_X\beta) = (-1)^k\langle\beta,\,D_{\flat(X)}\sharp_\omega\alpha\rangle 
 + \langle\alpha,\,D_{\flat(X)}\sharp_\omega\beta\rangle. 
\]
From Lemma \ref{sec5:lemma3} it follows that 
$\langle\alpha,\,D_{\flat(X)}\sharp_\omega\beta\rangle + (-1)^k\langle\beta,\,D_{\flat(X)}\sharp_\omega\alpha\rangle =\nabla_X\bigl(\Pi_k(\alpha,\beta)\bigr)$.
\qquad\qquad\qquad $\Box$
\vspace{0.5cm}

By using Proposition \ref{sec5:prop_Poisson connection} and the condition that $\Pi_{2m}(\mu,\,\mu)=\langle \mu,\,\hat{\mu}\rangle=1$, 
\begin{align*}
\langle\nabla_X\mu,\,\hat{\mu}\rangle &= 
\Pi_{2m}(\nabla_X\mu,\,\mu)=\frac{1}{2}\ \{\Pi_{2m}(\nabla_X\mu,\,\mu) + \Pi_{2m}(\mu,\,\nabla_X\mu)\}\\
&=\frac{1}{2}\ X(\Pi_{2m}(\mu,\,\mu))=0. 
\end{align*}
Since $\hat{\mu}=\sharp_{2m}(\mu)\ne 0$, we have $\nabla_X\mu=0$. Therefore, it holds that 
\begin{equation}\label{sec5:eqn_Poisson connection1}
 \nabla_X\bigl(\Pi_k(\alpha,\beta)\bigr)\mu = \nabla_X\bigl(\Pi_k(\alpha,\beta)\mu\bigr) = \nabla_X(\alpha\wedge\star\beta) 
 = (\nabla_X\alpha)\wedge\star\beta + \alpha\wedge\nabla_X(\star\beta) 
\end{equation}
for any differential $k$-form $\alpha,\,\beta$ and any vector field $X$ on $P$. On the other hand, 
\begin{equation}\label{sec5:eqn_Poisson connection2}
 \bigl\{\Pi_k(\nabla_X\alpha,\,\beta) + \Pi_k(\alpha,\,\nabla_X\beta)\bigr\}\mu = (\nabla_X\alpha)\wedge\star\beta + \alpha\wedge\star(\nabla_X\beta). 
\end{equation}
From (\ref{sec5:eqn_Poisson connection1}), (\ref{sec5:eqn_Poisson connection2}) and Proposition \ref{sec5:prop_Poisson connection}, we have 
\begin{equation}\label{sec5:eqn_nabla and star}
 \nabla_X\star = \star\nabla_X.
\end{equation}

By Lemma \ref{sec2:lemma_Poisson coboundary} and Proposition \ref{sec5:prop_relation}, 
the curl operator $\mathfrak{C}_{\mu}$ can be written as 
\begin{align*}
 \mathfrak{C}_{\mu}^k(K) &= ({\mu}^\sharp\circ\flat_{2m-k+1}\circ\partial^{2m-k}_\Pi\circ\sharp_{2m-k}\circ{\mu}^\flat)(K)\\
 &= (\sharp_{k-1}\circ\star\circ\flat_{2m-k+1}\circ\partial^{2m-k}_\Pi\circ\sharp_{2m-k}\circ\star\circ\flat_k)(K)\qquad (K\in \mathfrak{X}^k(P)). 
\end{align*}
Let $\{e_1,\cdots e_{2m}\}$ be a local frame on some local chart $U$ and $\{\theta^1,\cdots,\theta^{2m}\}$ its dual coframe field. 
For $K\in\mathfrak{X}^k(P)$, we put $\kappa=\flat_k(K)$ for the sake of simplicity. 
By Theorem \ref{sec4:thm_Poisson-Lichnerowics} and (\ref{sec5:eqn_nabla and star}),  
\begin{align*}
 \mathfrak{C}_{\mu}^k|_U(K) &= -\sum_{i=1}^{2m}\ \bigl(\sharp_{k-1}\circ\star\circ\flat_{2m-k+1}\circ (e_i\wedge D_{\theta^i})\circ\sharp_{2m-k}\circ\star\circ\flat_k\bigr)(K)\\ 
 &= -\sum_{i=1}^{2m}\ \bigl(\sharp_{k-1}\circ\star\circ \flat(e_i)\wedge (\flat_{2m-k}\circ D_{\theta^i}\circ\sharp_{2m-k})\bigr)(\star\kappa)\\
 &= -\sum_{i=1}^{2m}\ \bigl(\sharp_{k-1}\circ\star\circ \flat(e_i)\wedge (\nabla_{\sharp\theta^i}\circ\star)\bigr)(\kappa)\\
 &=  -\sum_{i=1}^{2m}\ \sharp_{k-1}\circ\star\bigl(\flat(e_i)\wedge (\star\nabla_{\sharp\theta^i}\kappa)\bigr).
\end{align*}
It follows from Lemma \ref{sec5:lemma2} and Proposition \ref{sec5:prop_relation2} that 
\[
\star\bigl(\flat(e_i)\wedge (\star\nabla_{\sharp\theta^i}\kappa)\bigr)=(-1)^{2m-k}i\ (e_i)\star(\star\nabla_{\sharp\theta^i}\kappa) 
= (-1)^{2m-k}i\ (e_i)\nabla_{\sharp\theta^i}\kappa. 
\]
Accordingly, we have 
\begin{equation}\label{sec5:eqn_curl operator}
 \mathfrak{C}_{\mu}^k|_U(K) = -(-1)^k\sum_{i=1}^{2m}\ (\sharp_{k-1}\circ i\ (e_i)\nabla_{\sharp\theta^i}\circ\flat_k)(K). 
\end{equation}
In (\ref{sec5:eqn_curl operator}), put $\eta_i=(\nabla_{\sharp\theta^i}\circ\flat_k)K$ and then 
\begin{align*}
 \langle \xi,\,\mathfrak{C}_{\mu}^k|_U(K)\rangle &= -(-1)^k\sum_{i=1}^{2m}\ \langle\xi,\,\sharp_{k-1}i\ (e_i)\eta_i\rangle 
  = \sum_{i=1}^{2m}\ \langle i\ (e_i)\eta_i,\,\sharp_{k-1}\xi\rangle \\
 &= \sum_{i=1}^{2m}\ \eta_i(e_i,\,\sharp_{k-1}\xi). 
\end{align*}
for any $(k-1)$-form $\xi$ on $P$. 
Since $\flat$ is invertible, there exists a local coframe field $\{\vartheta^1,\cdots,\vartheta^{2m}\}$ such that $e_i=\sharp\vartheta^i~(i=1,2,\cdots,2m)$. 
Therefore, 
\begin{align*}
 \sum_{i=1}^{2m}\ \eta_i(e_i,\,\sharp_{k-1}\xi) &= \sum_{i=1}^{2m}\ \langle\eta_i,\,\sharp_k(\vartheta^i\wedge\xi)\rangle 
 = -\sum_{i=1}^{2m}\ \langle \xi\wedge\vartheta^i,\,\sharp_k\eta_i\rangle \\
 &= -\sum_{i=1}^{2m}\ \langle \xi,\,i\ (\vartheta^i)D_{\theta^i}K\rangle, 
\end{align*}
which leads us to the following result:

\begin{thm}\label{sec5:thm_the curl operator}
Let $(P,\omega)$ be a $2m$-dimensional symplectic manifold equipped with a Poisson connection $D$. 
If the torsion of $D$ is vanishing, then the $k$-th curl operator $\mathfrak{C}_{\mu}^k$ with respect to the volume form $\mu=\frac{1}{m!}\wedge^m\omega$ 
on $P$ can be written in the form 
\[
 \mathfrak{C}_{\mu}^k|_U = -\sum_{i=1}^{2m}\ i\ (\flat(e_i))D_{\theta^i}. 
\]
on some local chart $U$. 
\end{thm}

It can be shown that 
Theorem \ref{sec5:thm_the curl operator} is independent of the choice of the local frame field 
by the way similar to Theorem \ref{sec4:thm_Poisson-Lichnerowics}. 
According to Darboux theorem, there exists a local neighborhood of each point of $P$ with coordinates $(x_1,\cdots,x_m, y_1,\cdots,y_m)$ 
on which the symplectic form $\omega$ is written in the form 
\[
\omega = \sum_{i=1}^m\ dx_i\wedge dy_i. 
\]
In Theorem \ref{sec5:thm_the curl operator}, if 
\[
 e_i = \begin{cases} \frac{\partial}{\partial x_i}\quad (i=1,\cdots,m),\\
                     \frac{\partial}{\partial y_i}\quad (i=m+1,\cdots,2m), \end{cases}
\]
and 
\[
 \theta^i = \begin{cases} dx_i\quad (i=1,\cdots,m),\\ 
              dy_i\quad (i=m+1,\cdots,2m), \end{cases}
\]
then $\mathfrak{C}_{\mu}^k$ is represented locally as follows:
\begin{cor}\label{sec5:cor_curl}
$\mathfrak{C}_{\mu}^k$ is given locally by 
\[
 \mathfrak{C}_{\mu}^k|_W = \sum_{i=1}^m\ i(dx_i)D_{dy_i} - i(dy_i)D_{dx_i} 
\]
with respect to Darboux coordinates $(W:x_1,\cdots,x_m,y_1,\cdots,y_m)$. 
\end{cor}

\begin{ex}
Let us consider a 2-dimensional symplectic manifold $\mathbb{R}^2$ with a symplectic form $\omega = \phi\,dx\wedge dy$, where $\phi$ is a smooth function 
$\phi(x,y)\ne 0$ on $\mathbb{R}^2$~(see Example \ref{sec3:ex_symplectic}). 
Let $D$ be a Poisson connection whose torsion is vanishing. 
For any vector field $X=f\frac{\partial}{\partial x} + g\frac{\partial}{\partial y}$ on $\mathbb{R}^2$, we have 
\begin{align*}
 -i\ (\flat(\partial/\partial x))D_{dx}(X) - i\ (\flat(\partial/\partial y))D_{dy}(X) &= -\phi\ (dy)(D_{dx}X) + \phi\ (dx)(D_{dy}X)\\
&= \frac{\partial f}{\partial x} + \frac{\partial g}{\partial y} + \phi\ (D_{dx}dy-D_{dy}dx)(X)\\
&= \mathrm{div}\,X + \phi\{dx,\,dy\}_\Pi(X)\\
&= \mathrm{div}\,X - \phi\,X\biggl(\frac{1}{\phi}\biggr).
\end{align*}
Here, we remark that $\{dx,\,dy\}_\Pi = -d\{x,\,y\} = -d\phi$. 
Therefore, it holds that 
\[
\mathfrak{C}_{\mu}(X) = -i\ (\flat(\partial/\partial x))D_{dx}(X) - i\ (\flat(\partial/\partial y))D_{dy}(X).
\]
In particular, when $\phi\equiv 1$, i.e., $\omega$ is the canonical symplectic form, 
\[
 i(dx)D_{dy}(X) - i(dy)D_{dx}(X) = -dx(D_{dx}X) + dy(D_{dy}X) = \mathfrak{C}_{\mu}(X), 
\]
which indicates Corollary \ref{sec5:cor_curl}. 
\end{ex}

From Proposition \ref{sec5:prop_relation} and \ref{sec5:prop_Brylinski},  
\[
  \mathfrak{C}_{\mu}^k = {\mu}^\sharp\circ d\circ{\mu}^\flat = \sharp_{k-1}\circ(\star\circ d\circ \star)\circ \flat_k 
               = (-1)^{k+1}\sharp_{k-1}\circ\delta_\Pi^k\circ\flat_k. 
\]
Using Theorem \ref{sec5:thm_the curl operator}, we can also get the following result: 
\begin{cor}
The differential operator $\delta_\Pi^\bullet$ is represented locally by 
\[
 \delta^k_\Pi|_U = \sum_{i=1}^{2m}\ i\ (e_i)\nabla_{\sharp\theta^i}.
\]
in terms of the symplectic connection $\nabla$ by (\ref{sec5:eqn_symplectic connection}).
\end{cor}
\subsection{Applications}

We let $(P,\,\omega)$ be a symplectic manifold of $2m$-dimension equipped with a Poisson connection $D$ whose curvature 2-section is $R_D$ and 
suppose that the torsion of $D$ is vanishing. As discussed in 5.2, $D$ induces a symplectic connection $\nabla$. 
The following lemma says that $D_\alpha$ for any 1-form $\alpha$ 
commutes with the map $\flat$ induced from $\omega$.
\begin{lem}\label{sec5:lem_flat}
For any $\alpha\in \Omega^1(P)$ and $X\in \mathfrak{X}^1(P)$, $D_\alpha(\flat(X)) = \flat(D_\alpha X)$. 
\end{lem}
{\em Proof}.~ Let $V$ be any vector field on $P$. On one hand, $\langle D_\alpha(\flat(X)),\,V\rangle = \langle \flat\circ\nabla_{\sharp\alpha}X,\,V\rangle 
 = \omega(\nabla_{\sharp\alpha}X,\,V)$. On the other hand, 
\begin{align*}
 \langle \flat(D_\alpha X),\,V\rangle &= \langle \flat(\sharp\nabla_{\sharp\alpha}\flat(X)),\,V\rangle = \langle \nabla_{\sharp\alpha}\flat(X),\,V\rangle \\
 &= (\sharp\alpha)\bigl(\omega(X,\, V)\bigr) - \omega\bigl(X,\,\nabla_{\sharp\alpha}V\bigr). 
\end{align*}
Since $\nabla$ is a symplectic connection, we obtain 
\[
\langle D_\alpha\flat(X)-\flat(D_\alpha X),\,V\rangle = \omega\bigl(\nabla_{\sharp\alpha}X,\,V\bigr)-(\sharp\alpha)\bigl(\omega(X,\, V)\bigr)
+\omega\bigl(X,\,\nabla_{\sharp\alpha}V\bigr) = 0, 
\]
which shows the lemma. 
\qquad\qquad\qquad\qquad\qquad\qquad\qquad\qquad\qquad\qquad\qquad\qquad\qquad\qquad\qquad $\Box$ 
\vspace{0.5cm}

Let $\{e_1,\cdots, e_{2m}\}$ be a local frame field on some local chart $U$ in $P$ and $\{\theta^1,\cdots,\theta^{2m}\}$ its dual coframe field. 
For each $\theta^i,\,\theta^j~(i,j=1,\cdots,2m)$, $D_{\theta^i}\theta^j$ are differential 1-forms on $U$. 
So, there exist a family of functions $\varGamma^{ij}_k$ on $U$ and then $D_{\theta^i}\theta^j$ 
can be written in the form $D_{\theta^i}\theta^j = \sum_{k}\varGamma^{ij}_k\theta^k$. By (\ref{sec4:contravariant derivative}), 
vector fields $D_{\theta^i}e_j$ on $U$ can be expressed as $D_{\theta^i}e_j = -\sum_k\varGamma^{ik}_j e_k$. 
More generally, if a vector field $X$ is locally written as $X=\sum_i\xi_i e_i$, then $D_{\theta^i}X$ is given by 
\begin{equation}
 D_{\theta^i}X = \sum_k\bigl\{\mathcal{L}_{\sharp\theta^i}\xi_k - \sum_j(D_{\theta^i}\theta^k)(e_j)\bigr\}e_k.
\end{equation}

We take $\mu=\frac{1}{m!}\wedge^m\omega$ as a volume form on $P$. 
Using Theorem \ref{sec4:thm_Poisson-Lichnerowics} and \ref{sec5:thm_the curl operator}, we obtain local explicit description for both 
$\partial_\Pi^a A$ and $\mathfrak{C}_{\mu}^{a}A$ of $A\in \mathfrak{X}^a(P)$. 
Therefore, we have 
\begin{align*}
 (\partial_\Pi^{a-1}|_U\!\circ\mathfrak{C}_{\mu}^{a}|_U)(A) &= -\sum_{i=1}^{2m} e_i\wedge D_{\theta^i}\biggl(-\sum_{j=1}^{2m} i\ (\flat(e_j))D_{\theta^j}A\biggr) 
 = \sum_{i,j=1}^{2m} e_i\wedge D_{\theta^i}\bigl(i\ (\flat(e_j))D_{\theta^j}\bigr)A \notag \\
 &= \sum_{i,j=1}^{2m} e_i\wedge i\ (D_{\theta^i}(\flat(e_j)))D_{\theta^j} A + \sum_{i,j=1}^{2m} e_i\wedge i\ (\flat(e_j))D_{\theta^i}D_{\theta^j} A 
\end{align*}
We remark that Proposition \ref{sec4:prop_inner product} is used in the last equality. 
On the other hand, by Proposition \ref{sec2:prop_inner product} and (\ref{sec4:eqn_distribution low}), 
\begin{align*}
 (\mathfrak{C}_{\mu}^{a+1}|_U\!\circ\partial_\Pi^{a}|_U)(A) &= \sum_{j=1}^{2m}i\ (\flat(e_j))D_{\theta^j}\biggl(\sum_{i=1}^{2m}e_i\wedge D_{\theta^i}A\biggr)\notag \\
 &= \sum_{i,j=1}^{2m}i\ (\flat(e_j))(D_{\theta^j}e_i\wedge D_{\theta^i} + e_i\wedge D_{\theta^j}D_{\theta^i})A\notag \\
 &= \sum_{i,j=1}^{2m}i\ (\flat(e_j))\biggl(-\sum_{a}\varGamma_i^{jk}e_a\wedge D_{\theta^i}A\biggr) 
                  + \sum_{i,j=1}^{2m}i\ (\flat(e_j))(e_i\wedge D_{\theta^j}D_{\theta^i})A\notag \\
 &= -\sum_{i,j=1}^{2m}i\ (\flat(e_j))\bigl(e_i\wedge D_{D_{\theta^j}\theta^i}A\bigr) + 
       \sum_{i,j=1}^{2m}\bigl\{e_i\wedge i\ (\flat(e_j))D_{\theta^j}D_{\theta^i}A + (-1)^k\omega(e_j,e_i)D_{\theta^j}D_{\theta^i}A\bigr\}\notag \\
 &= \sum_{i,j=1}^{2m}e_i\wedge i\ (\flat(e_j))\bigl(D_{\theta^j}D_{\theta^i}-D_{D_{\theta^j}\theta^i}\bigr)A 
      + (-1)^a \sum_{i,j=1}^{2m}\omega(e_i,e_j)\bigl(D_{\theta^i}D_{\theta^j}-D_{D_{\theta^i}\theta^j}\bigr)A. 
\end{align*}
Therefore, 
\begin{align*}
\Lambda_{\mu}^a|_UA &=(\mathfrak{C}_{\mu}^{a+1}\!\circ\partial_\Pi^{a}-\partial_\Pi^{a-1}\!\circ\mathfrak{C}_{\mu}^{a})A \\
&= -\sum_{i,j=1}^{2m}e_i\wedge i\ (\flat(e_j))\bigl(D_{\theta^i}D_{\theta^j}-D_{\theta^j}D_{\theta^i}+D_{D_{\theta^j}\theta^i}\bigr) A 
   - \sum_{i,j=1}^{2m} e_i\wedge i\ (D_{\theta^i}(\flat(e_j)))D_{\theta^j} A\\
  &\qquad + (-1)^a \sum_{i,j=1}^{2m}\omega(e_i,e_j)\bigl(D_{\theta^i}D_{\theta^j}-D_{D_{\theta^i}\theta^j}\bigr)A. 
\end{align*} 
Since it holds that $D_{\theta^j}\theta^i = D_{\theta^i}\theta^j - [\theta^i,\,\theta^j]_\Pi$ from the assumption, 
\begin{align*}
\Lambda_{\mu}^a|_UA &= -\sum_{i,j=1}^{2m}e_i\wedge i\ (\flat(e_j))\ R_D(\theta^i,\,\theta^j) A 
                            -  \sum_{i,j=1}^{2m} e_i\wedge i\ (D_{\theta^i}(\flat(e_j)))D_{\theta^j} A\\
 &\quad  + (-1)^a \sum_{i,j=1}^{2m}\omega(e_i,e_j)\bigl(D_{\theta^i}D_{\theta^j}-D_{D_{\theta^i}\theta^j}\bigr)A
  - \sum_{i,j=1}^{2m}e_i\wedge i\ (\flat(e_j))D_{D_{\theta^i}\theta^j}A\\
 &= -\sum_{i,j=1}^{2m}e_i\wedge i\ (\flat(e_j))\ R_D(\theta^i,\,\theta^j) A - \sum_{i,j=1}^{2m} e_i\wedge i\ \bigl(D_{\theta^i}(\flat(e_j))\bigr)D_{\theta^j}A\\                      
 &\quad + (-1)^a \sum_{i,j=1}^{2m}\omega(e_i,e_j)\bigl(D_{\theta^i}D_{\theta^j}-D_{D_{\theta^i}\theta^j}\bigr)A
          - \sum_{i,j=1}^{2m}e_i\wedge i\ (\flat(e_j))D_{\sum_k\varGamma^{ij}_k \theta^k}A \\ 
 &= -\sum_{i,j=1}^{2m}e_i\wedge i\ (\flat(e_j))\ R_D(\theta^i,\,\theta^j) A 
                            - \sum_{i,j=1}^{2m} e_i\wedge i\ \bigl(D_{\theta^i}(\flat(e_j))\bigr)D_{\theta^j} A \\
   &\quad + (-1)^a \sum_{i,j=1}^{2m}\omega(e_i,e_j)\bigl(D_{\theta^i}D_{\theta^j}-D_{D_{\theta^i}\theta^j}\bigr)A 
                            + \sum_{i,j=1}^{2m}e_i\wedge i\ \bigl((\flat(D_{\theta^i}e_j))\bigr)D_{\theta^j}A. 
\end{align*} 
Using Lemma \ref{sec5:lem_flat} in the last equality, we have 
\[
 \Lambda_{\mu}^a|_UA =(-1)^a\sum_{i,j=1}^{2m}\omega(e_i,\,e_j)(D_{\theta^i}D_{\theta^j} - D_{D_{\theta^i}\theta^j})A 
    - \sum_{i,j=1}^{2m}e_i\wedge i\ (\flat(e_j))\ R_D(\theta^i,\,\theta^j)A. 
 \]
Put $D_{ij}=D_{\theta^i}D_{\theta^j} - D_{D_{\theta^i}\theta^j}$ for each $i,j=1,2,\cdots, 2m$ and then 
\begin{align*}
 \sum_{i,j=1}^{2m}\omega(e_i,\,e_j)\bigl(D_{\theta^i}D_{\theta^j}-D_{D_{\theta^i}\theta^j}\bigr) 
 &= \sum_{i<j}^{2m}\omega(e_i,\,e_j)(D_{ij}-D_{ji})\\ 
 &= \sum_{i<j}^{2m}\omega(e_i,\,e_j)(D_{\theta^i}D_{\theta^j}-D_{\theta^j}D_{\theta^i}-D_{[\theta^i, \theta^j]_\Pi})\\
 &= \sum_{i<j}^{2m}\omega(e_i,\,e_j)R_D(\theta^i,\,\theta^j).
\end{align*}
Consequently, we find that the modular operator with respect to $\mu=\frac{1}{m!}\wedge^m\omega$ is locally expressed in the form 
 \begin{equation*}\label{sec5:eqn_main3}
  \Lambda_{\mu}^a|_UA = (-1)^a\sum_{i<j}^{2m}\omega(e_i,\,e_j)\ R_D(\theta^i,\,\theta^j)A - \sum_{i,j=1}^{2m}e_i\wedge i\ (\flat(e_j))\ R_D(\theta^i,\,\theta^j)A 
 \end{equation*}
By the fact the modular vector field $\Xi_{\mu}$ with respect to $\mu$ is zero, we can get the following result: 

\begin{thm}\label{sec5:thm_main} 
Let $\{e_1,\cdots, e_{2m}\}$ be a local frame field on some local chart $U$ in $P$ and $\{\theta^1,\cdots,\theta^{2m}\}$ its dual coframe field. 
If $P$ has a Poisson connection $D$ whose torsion is vanishing, then $R_D$ satisfies the following equation:
 \[
 (-1)^a\sum_{i<j}^{2m}\omega(e_i,\,e_j)\ R_D(\theta^i,\,\theta^j)\ A 
    = \sum_{i,j=1}^{2m}e_i\wedge i\ (\flat(e_j))\ R_D(\theta^i,\,\theta^j)\ A. 
 \]
for any $A\in\mathfrak{X}^a(P)$. 
\end{thm}
Especially, if we choose Darboux coordinates $(x_1,\cdots,x_m,y_1,\cdots,y_m)$ as $U$, the above formula is given by 
\begin{align*}
 \sum_{i=1}^m R_D(dx_i, dy_i)\ A &= \sum_{i=1}^m \Bigl\{\sum_{j=1}^m i(dx_j)R_D(dx_i, dy_j)
            -\sum_{j\ne i} i(dy_j)R_D(dx_i, dx_j)\Bigr\}\ A\wedge \frac{\partial}{\partial x_i}\\
 &\qquad + \sum_{i=1}^m \Bigl\{\sum_{j=1}^m i(dy_j)R_D(dx_j, dy_i)
            - \sum_{j\ne i} i(dx_j)R_D(dy_j, dy_i)\Bigr\}\ A\wedge \frac{\partial}{\partial y_i}. 
\end{align*}
\vspace{0.8cm}

Next, let us observe how the modular operator for a symplectic manifold $(P,\,\omega)$ can be locally described when we change the volume form on $P$. 
If $f$ is a non-vanishing function, then it holds that, for any $A\in\mathfrak{X}^a(P)$ 
\begin{equation*}
 [A,\,\log|f|]_{\mathrm{S}} = (-1)^{a-1}i\ (d\log|f|)A 
\end{equation*}
by using (\ref{sec2:eqn_Schouten}). From this, the curl operator with respect to the volume form $\nu = f\mu$ is calculated as 
\begin{align*}
 \mathfrak{C}_{\nu}^a A &= \nu^{\sharp}(di(A)\nu) = \nu^{\sharp}(fi(A)\mu) + \nu^{\sharp}(df\wedge i(A)\mu) 
              = i(di(A)\mu)\ \hat{\mu} + i(d\log|f|\wedge i(A)\mu)\ \hat{\mu}\\
 &= \mathfrak{C}_{\mu}^a A + i(d\log|f|)\ \mu^{\sharp}(i(A)\mu) = \mathfrak{C}_{\mu}^a A + i(d\log|f|)\ A\\ 
 &= \mathfrak{C}_{\mu}^a A + (-1)^{a-1}[A,\,\log|f|]_{\mathrm{S}}.
\end{align*}
Applying (\ref{sec3:eqn_Schouten and curl}) and Theorem \ref{sec5:thm_the curl operator} to the second term in the last equality, we have 
\begin{equation}\label{sec5:eqn_volume form changed}
 \mathfrak{C}_{\nu}^a|_U A = \mathfrak{C}_{\mu}^a|_U A - \sum_{i=1}^{2m}(\mathcal{L}_{\sharp\theta^i}\log|f|)\ i(\flat(e_i))\ A 
\end{equation}
on a local chart $U$ of $P$. 
\begin{thm}\label{sec5:thm_main2} 
Let $\{e_1,\cdots, e_{2m}\}$ be a local frame field on some local chart $U$ in $P$ and $\{\theta^1,\cdots,\theta^{2m}\}$ its dual coframe field. 
If $P$ has a Poisson connection $D$ whose torsion is vanishing, then the modular operator with respect to the volume form $\nu=f\mu$ can be locally described 
in the explicit formula 
\begin{align*}
\Lambda_{\nu}^a|_UA &= (-1)^a\sum_{i,j=1}^{2m}(\mathcal{L}_{\sharp\theta^i}\log|f|)\ \omega(e_i,e_j)D_{\theta^j}A \\ 
 &\qquad + \sum_{i,j=1}^{2m}\biggl\{\sum_{k=1}^{2m}(\mathcal{L}_{\sharp\theta^k}\log|f|)\varGamma_k^{i,j}
  - \mathcal{L}_{\sharp\theta^i}\mathcal{L}_{\sharp\theta^j}\log|f|\biggr\}e_i\wedge i(\flat(e_j))A. 
\end{align*}
\end{thm}

\noindent {\em Proof.}~ By (\ref{sec5:eqn_volume form changed}), we have 
\begin{align*}
 (\mathfrak{C}_{\nu}|_U\circ\partial_\Pi|_U)\ A &= (\mathfrak{C}_{\mu}|_U\circ\partial_\Pi|_U)\ A 
               - \sum_{i=1}^{2m}(\mathcal{L}_{\sharp\theta^i}\log|f|)\ i(\flat(e_i))\ \partial_\Pi|_U A, \\
 (\partial_\Pi|_U\circ\mathfrak{C}_{\nu}|_U)\ A &= (\partial_\Pi|_U\circ\mathfrak{C}_{\mu}|_U)\ A 
               - \sum_{i=1}^{2m}\partial_\Pi|_U\Bigl((\mathcal{L}_{\sharp\theta^i}\log|f|)\ i(\flat(e_i))\ A\Bigr). 
\end{align*}
Therefore, 
\begin{equation}\label{sec5:eqn_main2}
\Lambda_{\nu}^a|_U A = -\sum_{i=1}^{2m}\Bigl\{F_i\ i(\flat(e_i))\partial_\Pi A - \partial_\Pi(F_i\flat(e_i)A)\Bigr\},
\end{equation}
where $F_i=\mathcal{L}_{\sharp\theta^i}\log|f|~(i=1,\cdots,2m)$. We remark that $\Lambda_{\mu}|_U=\mathfrak{C}_{\mu}|_U\circ\partial_\Pi|_U 
 - \partial_\Pi|_U\circ\mathfrak{C}_{\mu}|_U = 0$. 
By Theorem \ref{sec4:thm_Poisson-Lichnerowics}, 
\[
F_i\ i(\flat(e_i))\partial_\Pi|_U A = -F_i\sum_{j=1}^{2m}\{e_j\wedge i(\flat(e_i))D_{\theta_j}A + (-1)^a\omega(e_i,\,e_j)D_{\theta^j}A\} 
\]
and 
\[
\partial_\Pi|_U(F_i\flat(e_i)A) = -\sum_{j=1}^{2m}e_j\wedge\{F_i D_{\theta^j}i(\flat(e_i))A + (\mathcal{L}_{\sharp\theta^j}F_i)i(\flat(e_i))A\}. 
\]
Using Proposition \ref{sec4:prop_inner product} and Lemma \ref{sec5:lem_flat}, we compute $\Lambda_{\nu}^a|_U A$ as 
\begin{align*}
(\ref{sec5:eqn_main2}) &= \sum_{i,j=1}^{2m}F_ie_j\wedge i\bigl(\flat(D_{\theta_j}e_i)\bigr)\ A 
      - \sum_{i,j=1}^{2m}(\mathcal{L}_{\sharp\theta^j}F_i)e_j\wedge i(\flat(e_i))\ A + (-1)^a\sum_{i,j=1}^{2m}F_i \omega(e_i,e_j)D_{\theta^j}\ A\\
 &= \sum_{i,j=1}^{2m}F_i\sum_{k=1}^{2m}(D_{\theta^j}\theta^k)(e_i)e_j\wedge i(\flat(e_k))\ A 
      - \sum_{i,j=1}^{2m}(\mathcal{L}_{\sharp\theta^j}F_i)e_j\wedge i(\flat(e_i))\ A + (-1)^a\sum_{i,j=1}^{2m}F_i \omega(e_i,e_j)D_{\theta^j}\ A. 
\end{align*}
Consequently, Theorem \ref{sec5:thm_main2} is shown. \qquad\qquad\qquad\qquad\qquad\qquad\qquad\qquad\qquad\qquad\qquad\qquad $\Box$
\vspace{0.5cm}

By (\ref{sec3:eqn_modular operator}), we can specify locally the Schouten bracket of the modular vector field $\Xi_{\nu}$ with respect to $\nu$ 
and $A\in\mathfrak{X}^a(P)$ as 
\begin{align*}
 [\Xi_{\nu},\,A]_{\rm S}|_U &= \sum_{i,j=1}^{2m}(\mathcal{L}_{\sharp\theta^i}\log|f|)\ \omega(e_i,e_j)D_{\theta^j}A \\ 
 &\qquad + (-1)^a\sum_{i,j=1}^{2m}\biggl\{\sum_{k=1}^{2m}(\mathcal{L}_{\sharp\theta^k}\log|f|)\varGamma_k^{i,j}
  - \mathcal{L}_{\sharp\theta^i}\mathcal{L}_{\sharp\theta^j}\log|f|\biggr\}e_i\wedge i(\flat(e_j))A. 
\end{align*}
in terms of Poisson connection whose torsion is vanishing. 
\begin{ex}
 We consider a symplectic manifold $\mathbb{R}^2$ equipped with a symplectic form $\omega=(1/\phi)\ dx\wedge dy$ 
 and compute the modular operator $\Xi_{\nu}X$ with respect to the volume form $\nu=dx\wedge dy$ of $X=f\frac{\partial}{\partial x}+g\frac{\partial}{\partial y}$
 ~(see Example \ref{sec3:ex_2-dimensional case}). Remark that $\nu$ is obtained from $\omega$ multiplied by $\phi$. 
 Let $D$ be a Poisson connection whose torsion is vanishing. We choose $\{\frac{\partial}{\partial x},\,\frac{\partial}{\partial y}\}$ 
 as a local frame field and put $e_1=\partial_x=\frac{\partial}{\partial x},\,e_2=\partial_y=\frac{\partial}{\partial y}$. 
 Then, by a direct computation we have  
\begin{align*}
&-\sum_{i,j=1}^{2m}(\mathcal{L}_{\sharp\theta^i}\mathcal{L}_{\sharp\theta^j}\log|f|)\ e_i\wedge i(\flat(e_j)) = 
\biggl(-\frac{\partial^2\phi}{\partial x\partial y}f - \frac{\partial^2\phi}{\partial y^2}g\biggr)\ \frac{\partial}{\partial x} + 
\biggl(\frac{\partial^2\phi}{\partial x^2}f + \frac{\partial^2\phi}{\partial x\partial y}g\biggr)\ \frac{\partial}{\partial y},\\ 
&\sum_{i,j=1}^{2m}\sum_{k=1}^{2m}(\mathcal{L}_{\sharp\theta^k}\log|f|)(D_{\theta^i}\theta^j)(e_k)\ e_i\wedge i(\flat(e_j))\ A \\ 
&\qquad\qquad= \frac{1}{\phi}\biggl\{-g\frac{\partial\phi}{\partial y}(D_{dx}dx)(\partial_x) + g\frac{\partial\phi}{\partial y}(D_{dx}dx)(\partial_y) 
+ f\frac{\partial\phi}{\partial y}(D_{dx}dy)(\partial_x)-f\frac{\partial\phi}{\partial x}(D_{dx}dy)(\partial_y)\biggr\}\ \frac{\partial}{\partial x}\\
&\qquad\qquad+ \frac{1}{\phi}\biggl\{-g\frac{\partial\phi}{\partial y}(D_{dy}dx)(\partial_x) + g\frac{\partial\phi}{\partial y}(D_{dy}dx)(\partial_y) 
+ f\frac{\partial\phi}{\partial y}(D_{dy}dy)(\partial_x)-f\frac{\partial\phi}{\partial x}(D_{dy}dy)(\partial_y)\biggr\}\ \frac{\partial}{\partial y}
\end{align*}
and 
\begin{align*}
&-\sum_{i,j=1}^{2m}(\mathcal{L}_{\sharp\theta^i}\log|f|) \omega(e_i,e_j)D_{\theta^j}A \\
=&-\frac{1}{\phi}\biggl\{f\frac{\partial\phi}{\partial x}(D_{dx}dx)(\partial_x) + g\frac{\partial\phi}{\partial x}(D_{dx}dx)(\partial_y) 
+ f\frac{\partial\phi}{\partial y}(D_{dy}dx)(\partial_x)+g\frac{\partial\phi}{\partial y}(D_{dy}dx)(\partial_y)\biggr\}\ \frac{\partial}{\partial x}\\
&-\frac{1}{\phi}\biggl\{f\frac{\partial\phi}{\partial x}(D_{dx}dy)(\partial_x) + g\frac{\partial\phi}{\partial x}(D_{dx}dy)(\partial_y) 
+ f\frac{\partial\phi}{\partial y}(D_{dy}dy)(\partial_x)+g\frac{\partial\phi}{\partial y}(D_{dy}dy)(\partial_y)\biggr\}\ \frac{\partial}{\partial y}\\
& + \biggl(\frac{\partial f}{\partial x}\frac{\partial\phi}{\partial y}-\frac{\partial f}{\partial y}\frac{\partial\phi}{\partial x}\biggr)\ \frac{\partial}{\partial x}
 + \biggl(\frac{\partial g}{\partial x}\frac{\partial\phi}{\partial y}-\frac{\partial g}{\partial y}\frac{\partial\phi}{\partial x}\biggr)\ \frac{\partial}{\partial y}. 
\end{align*}
From the assumption that $D$ is a Poisson connection it follows that 
\begin{equation}\label{sec5:eqn2_main2}
 (D_{dx}dx)(\partial_x) + (D_{dx}dy)(\partial_y) = -\frac{\partial\phi}{\partial y},\quad 
 (D_{dy}dx)(\partial_x) + (D_{dy}dy)(\partial_y) =  \frac{\partial\phi}{\partial x}. 
\end{equation}
By (\ref{sec5:eqn2_main2}) and the fact that the torsion of $D$ is vanishing, the right-hand side of the formula in Theorem \ref{sec5:thm_main2} is calculated to be 
\begin{equation*}
 \biggl(\frac{\partial f}{\partial x}\frac{\partial\phi}{\partial y}-\frac{\partial f}{\partial y}\frac{\partial\phi}{\partial x}
 -f\frac{\partial^2\phi}{\partial x\partial y}-g\frac{\partial^2\phi}{\partial y^2}\biggr)\ \frac{\partial}{\partial x}
 + \biggl(\frac{\partial g}{\partial x}\frac{\partial\phi}{\partial y}-\frac{\partial g}{\partial y}\frac{\partial\phi}{\partial x}
 + f\frac{\partial^2\phi}{\partial x^2}+g\frac{\partial^2\phi}{\partial x\partial y}\biggr)\ \frac{\partial}{\partial y}, 
\end{equation*}
which certainly gives us the same result as Example \ref{sec3:ex_2-dimensional case}. 
\end{ex}
\vspace{0.8cm}

\noindent{\bf Acknowledgments.}~The author would like to express his deepest 
gratitude to Emeritus Professor Toshiaki Kori, Professor Hiroaki Yoshimura of Waseda University and Noriaki Ikeda of Ritsumeikan University for 
useful discussion and helpful comments. 
He also wishes to thank Waseda University for the hospitality while part of the work was being done. 

\noindent 
{\small 
{\em Yuji HIROTA}\\
{Azabu University}\\
\textit{Email} : \texttt{hirota@azabu-u.ac.jp}
}

\begin{thebibliography}{99}
 \bibitem{Bpan03} M. Berger: {\em A panoramic view of Riemannian geometry}, Springer-Verlag, Berlin, 2003. 

 \bibitem{BCGRSsym06} P. Bieliavsky, M. Cahen, S. Gutt, J. Rawnsley and L. Schwarchh$\ddot{\rm o}$fer: {\em Symplectic connections}, 
  Int. J. Geom. Methods Mod. Phys. 3 (2006), no. 3, 375--420. 
 
 \bibitem{Bdiff88} J. -L. Brylinski: {\em A differential complex for Poisson manifolds}, J. Differential Geom. 28 (1988), no. 1, 93--114. 
 
 \bibitem{CWgeo99} A. Cannas da Silva and A. Weinstein: 
 {\em Geometric models for noncommutative algebras}, American Mathematical Society, 
 Providence, RI, 1999. 

 \bibitem{DZpoi05} J.-P. Dufour and N. T. Zung: {\em Poisson Structures and Their Normal Forms}. 
  Prog. in Math. {242}, Birkh$\ddot{\mathrm{a}}$user, Basel, 2005. 

 \bibitem{Fcon00} R. L. Fernandes: {\em Connections in Poisson geometry I: Holonomy and invariants}, J. Differential Geom. 54 (2000), no. 2, 303--365. 

 \bibitem{Acrochet85} J.-L. Koszul: {\em Crochet de Schouten-Nijenhuis et cohomologie}, The mathematical heritage of $\acute{\rm E}$lie Cartan (Lyon, 1984), 
 As$\acute{\rm e}$rique 1985, Num$\acute{\rm e}$ro Hors S$\acute{\rm e}$rie, 257--271. 

 \bibitem{LPVpoi12} C. Laurent-Gengoux, A. Pichereau, P. Vanhaecke: {\em Poisson strucutures}, Grundlehren der Mathematischen Wissenschaften, 347, 
 Springer, Heidelberg, 2013. 

 \bibitem{Jrie02} J. Jost: {\em Riemannian geometry and geometric analysis}, Third edition. Universitext. Springer-Verlag, Berlin, 2002. 

\end{thebibliography}
\end{document}